\documentclass[12pt]{article}
\usepackage{preprint-layout}
\usepackage{preprint-notation}

\title{Hybridized CutFEM for\\ Elliptic Interface Problems}
\author{Erik Burman, Daniel Elfverson, Peter Hansbo,\\
Mats G. Larson and
Karl Larsson
}
\date{}

\begin{document}

\maketitle
\begin{abstract}
We design and analyze a hybridized cut finite element method for elliptic interface problems. In this method very general meshes can be coupled over internal unfitted interfaces, through a skeletal variable, using a Nitsche type approach. We discuss how optimal error estimates for the method are obtained using the tools of cut finite element methods and prove a condition number estimate for the Schur complement. Finally, we present illustrating numerical examples.
\end{abstract}

\section{Introduction}
The solution of heterogeneous problems, for instance problems where some
physical parameter has important variation within the computational
domain, remains a challenging problem in computational science. Indeed
special care must be taken in the design of methods to ensure that
efficiency and accuracy do not degenerate in the presence of high contrast 
inclusions. An additional layer of complexity is
added if the internal structure, e.g. the position of the inclusion, is
one of the problem unknowns and domain boundaries or internal
interfaces must be modified during the computation. In this situation
it is an advantage if remeshing of the domain can be avoided, while
the equations still are solved accurately on a variety of
configurations \cite{VM17, BEHLL18}. When one interface separates two
computational domains and the problem size is moderate it is
reasonable to use a monolithic solution strategy. However, as the number of inclusions
increases it becomes advantageous to resort to domain decomposition so that the
problem can be solved in parallel. 

Recently there has been a surge of activity in the design and analysis
of hybridized methods, that is nonconforming methods that have different
approximating polynomials
defined in the bulk of the elements and on the
skeleton of the computational mesh. The skeleton variable plays the
role of a mortar variable, either for a Neumann or a Dirichlet trace. 
Typically the interior degrees of
freedom of each cell can easily be eliminated using static
condensation, thereby reducing the size of the linear system. This is
particularly appealing for high order approximation methods where 
each volume element hosts a relatively large number of degrees of
freedom. An important feature of hybridized methods is that they allow for
fairly general element shapes and there exists an important literature
on methods defined on general polygonal/polyhedral meshes. Examples of
such methods include the Hybridized Discontinuous Galerkin Method
\cite{CoGoLa09}, the Virtual Element Method \cite{BBMR14} and the
Hybridized High Order method \cite{DPEL14}. In many cases these
methods are closely related (see \cite{CDPE16} and
references therein).

In this contribution we design a hybridized finite element method
within the cut finite element method (CutFEM) paradigm, see \cite{BuHa2012,BuClHaLaMa15}. This means 
that the computational mesh is independent of the geometry and internal 
interfaces, for example the computational mesh can
remain completely structured. If the underlying problem has some
special structure dividing the computational domain in subdomains, for instance
defined by grains with a specific permeability or microstructure, 
the domain boundaries are
allowed to cut through the background mesh in a close to arbitrary
fashion.
By including boundary or interface defining terms in the
variational formulation the method essentially eliminates the
meshing problem. This is particularly convenient in shape optimization
problems or inverse identification problems, where the position of the
interfaces change during an optimization process. The introduction of
a skeletal variable makes it possible to eliminate internal degrees of
freedom in a parallel fashion so that the linear system can be solved
by iterating on the Schur complement. Optimal stability and accuracy
is obtained using the tools developed within the CutFEM paradigm.
Below we will also show that the resulting Schur complement has a
condition number that is bounded independently of the mesh interface intersection.

\paragraph{Brief Review of Cut Finite Elements.}
CutFEMs were originally introduced by Hansbo and
Hansbo \cite{HH02} as an
alternative to extended finite element methods, using Nitsche's method
as mortaring over an unfitted interface. The ideas were extended to
composite meshes, 
Chimera-style, by Hansbo, Hansbo and Larson \cite{HHL03}
and then to fictitious domain problems by Burman and Hansbo \cite{BuHa2012}. In a
parallel development \cite{ORG09}, Olshanski, Reusken and Grande developed a cut
finite element method for the discretization of PDEs on surfaces using
the trace of the (discrete) surface on a finite element bulk mesh as
the computational mesh. Additional stability can be obtained by adding
suitable stabilization terms \cite{Bu10, BHL15,BHLM16,LaZa2017,GLR18} and 
the methods are then comparable to standard finite element methods on meshed 
geometries, both with respect to conditioning and accuracy. Further developments 
include techniques for the handling of curved interfaces or boundaries
\cite{LR17,BHL18}, PDEs on composite surfaces \cite{HJLL17} and transport 
problems in fractured mixed dimensional domains \cite{BHLL18}.
For an overview of the ideas
behind the CutFEM paradigm see \cite{BuClHaLaMa15} and for a
collection of essays giving a snapshot of the state of the art we
refer to \cite{BBLO18}.

\paragraph{New Contributions.}
We develop a hybridized CutFEM for an elliptic model 
problem with piecewise constant coefficients defined on a partition of the 
domain. The union of the boundaries of the subdomains is called the skeleton 
and hybridization consists of adding a solution field representing the 
solution on the skeleton. In the proposed method each subdomain of the bulk 
and each component of the skeleton is equipped with its own finite element 
mesh and space. The mesh may be constructed using a cut technique or a 
standard mesh. The hybridization leads to a convenient formulation which 
also naturally facilitates elimination of the bulk degrees of freedom using a 
Schur complement formulation. 

We develop a stability and error analysis where we in particular show that optimal 
order a priori error estimates holds for both and that the condition number of the 
Schur complement is $O(d^{-1} h^{-1} )$ where $d$ is the diameter of the 
subdomain. We cover very general choices of the subdomains and one instance 
of our method and analysis is a $p$ based discontinuous Galerkin method on polygonal 
elements.

\paragraph{Earlier Work.}

Hybridized methods are commonly used and many versions have been proposed, notably in the setting of discontinuous Galerkin methods; for an overview, see Cockburn et al. \cite{CoGoLa09}.
The particular version we employ here was first proposed using meshed subdomains by Egger \cite{Eg09}, who later extended it to incompressible flow \cite{EgWa12,EgWa13} and convection--diffusion problems \cite{ErSch10}; cf. also K\"onn\"o and Stenberg \cite{KoSt12} where it was employed for solving the Brinkman problem. Independently, the same approach has been proposed by Oikawa and Kikuchi \cite{OiKi10} and further developed by Miyashita and Saito \cite{MiSa17}.

 \paragraph{Outline.} In Section 2 we formulate the hybridized CutFEM. 
 In Section 3 we derive stability and optimal order error estimates.
 In Section 4 we eliminate the bulk fields using the Schur complement and we derive a bound for the condition number of the Schur complement.
 In Section 5 we present numerical results.
 
\section{The Hybridzed Cut Finite Element Method}
 
\subsection{The Model Problem}

\begin{figure}\centering
\begin{subfigure}[t]{0.35\linewidth}\centering
\includegraphics[width=\linewidth]{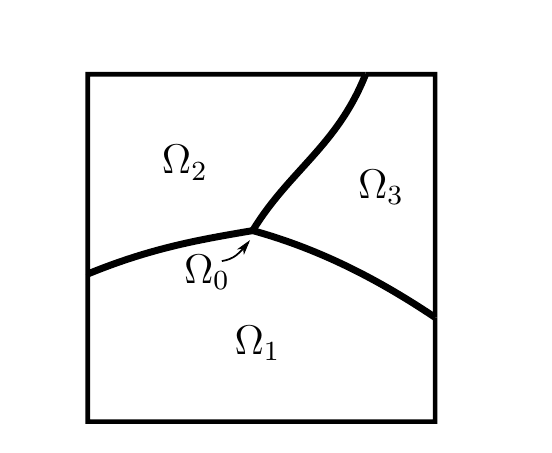}
\end{subfigure}
\begin{subfigure}[t]{0.35\linewidth}\centering
\includegraphics[width=\linewidth]{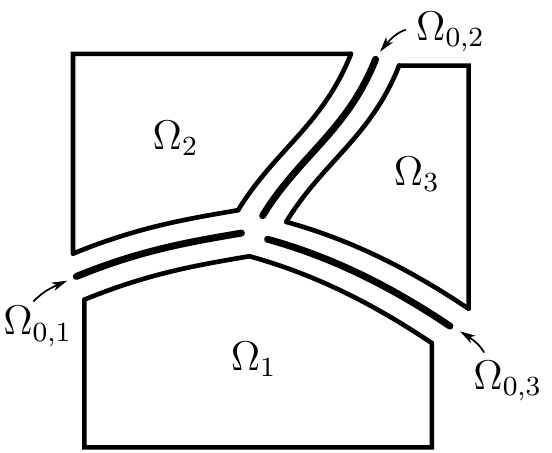}
\end{subfigure}
\caption{Schematic illustration of a partition of the domain $\Omega$ into a set of subdomains $\{\Omega_i\}_{i=1}^N$ and of the skeleton $\Omega_0$ with its partition into a set of skeleton subdomains $\{\Omega_{0,i}\}_{i=1}^{N_0}$.}
\label{fig:schematic-hybridized}
\end{figure}

\paragraph{Domain.} In hybridized methods the domain $\Omega$ is partitioned into a set of subdomains $\{\Omega_i\}_{i=1}^N$ and the interfaces between between the subdomains constitute a skeleton $\Omega_0$, see the illustration in Figure~\ref{fig:schematic-hybridized}.
Formally we assume the following:
\begin{itemize}
\item Let $\Omega \subset \IR^d$, with $d=2$ or $3$,  be a domain with piecewise smooth Lipschitz boundary $\partial \Omega$. Recall that $\partial \Omega$ is piecewise smooth if 
it is the union of a finite number of smooth $d-1$ dimensional manifolds with boundary and 
Lipschitz if there is a finite set of local coordinate systems in which the boundary can be described as a Lipschitz function.
%

\item Let $\{\Omega_i\}_{i=1}^N$ be a partition of 
$\Omega$ into $N$ subdomains $\Omega_i$ with piecewise smooth Lipschitz 
boundaries $\partial \Omega_i$.

\item Let $\Omega_0 = \cup_{i=1}^N \partial \Omega_i \setminus \partial \Omega$ be 
the skeleton of the partition. Note that there is a partition of $\Omega_0$ into 
smooth $d-1$ manifolds with boundary 
\begin{equation}
\Omega_0 = \cup_{i,j=1}^N \partial \Omega_i \cap \partial \Omega_j 
= \cup_{i=1}^{N_0} \Omega_{0,i}
\end{equation} 
where $N_0$ is the number of nonempty intersections $\partial \Omega_i \cap \partial \Omega_j$.
\end{itemize}

\paragraph{Model Problem.}
We consider the following hybridized formulation of the Poisson problem: 
find $u_0: \Omega_0 \rightarrow \IR$ and for $i=1,\dots,N$, 
$u_i:\Omega_i \rightarrow \IR$ such that 
\begin{alignat}{3}\label{eq:strong-a}
-\nabla \cdot a_i \nabla u_i &= f_i &\qquad &\text{in $\Omega_i$}
\\ \label{eq:strong-b}
\llbracket n\cdot a \nabla u \rrbracket &= 0 & \qquad &\text{on $\Omega_0$}
\\ \label{eq:strong-c}
[u]_i &=0  & \qquad &\text{on $\partial \Omega_i \cap \Omega_0$}
\\ \label{eq:strong-d}
u_i&=0 & \qquad  &\text{on $\partial \Omega_i \cap \partial \Omega$}
\end{alignat}
Here $a_i$, $i=0,\dots,N,$ are positive constants and the jumps operators are 
defined by
\begin{equation}
[v]|_{\partial \Omega_i \cap \Omega_0} = v_i - v_0, 
\qquad 
\llbracket n \cdot a \nabla v \rrbracket|_{\partial \Omega_i \cap \partial \Omega_j} 
= n_i \cdot a_i \nabla v_i +  n_j \cdot a_j \nabla v_j
\end{equation}
where $n_i$ is the exterior unit normal to $\Omega_i$. 

This problem is well posed, with exact solution $u \in H^1_0(\Omega)$
such that $u_i = u\vert_{\Omega_i}$, $i=0,\hdots,N$. Below we will
assume the additional regularity $u_i \in
H^{\frac32}(\Omega_i)$. We define the spaces
\begin{equation}\label{eq:reg-assumptions}
V_0 := \bigoplus_{i=1}^{N_0} L^2(\Omega_{0,i}),\quad V_{1,N} :=
\bigoplus_{i=1}^{N} H^{\frac32}(\Omega_i)
\end{equation}
and the global space
\begin{equation}
W := V_{0} \oplus V_{1,N}.
\end{equation}
We will use the notation $(\cdot,\cdot)_X$ for the $L^2$-scalar
product over $X \subset \mathbb{R}^s$ where either $s=d$ or
$s=d-1$. The associated norm will be denoted $\|\cdot\|_X$.
We will also use the following weighted $L^2$-norm, $\| v \|^2_{\omega, a} = (a v, v)_\omega$, for 
$a\in L^\infty(\omega)$ with $a > 0$.

\begin{rem} We note that solutions to certain problems of the form (\ref{eq:strong-a}-\ref{eq:strong-d}) 
may have regularity arbitrarily close to $H^1(\Omega)$, see \cite{K75} where a problem with intersecting interfaces and discontinuous coefficients are considered. Thus (\ref{eq:reg-assumptions}) does not hold in general. For smooth interfaces the assumption holds and of course for the situation when the coefficient is 
globally constant. Both of these cases are of practical importance since hybridization is used to conveniently 
enable elimination of the subdomain unknowns using a Schur complement procedure, see Section 3.
Deriving error bounds with lower regularity assumptions requires a more involved analysis 
and we refer to the techniques developed in \cite{BHL18b} for CutFEM approximations of boundary value problems with mixed Dirichlet/Neumann boundary conditions and minimal regularity requirements.
\end{rem}
\subsection{The Method}

\paragraph{Finite Element Spaces.}
On each subdomain $\Omega_i$, $1\leq i \leq N$, and on each skeleton subdomain 
$\Omega_{0,i}$, $1\leq i \leq N_0$, we construct a separate finite element space.
\begin{itemize}
\item Let $O \in \{ \Omega_{0,i}\}_{i=1}^{N_0} \cup \{\Omega_i\}_{i=1}^N$, i.e., 
$O$ is either one of the $d$ dimensional bulk subdomains $\Omega_i$ or one of 
the $d-1$ dimensional components of the skeleton $\Omega_0$.
\item For each $O$ there is a $d$ dimensional domain $U_O$ such that 
$O \subseteq U_O$. For each $h \in (0,h_{\max,i}]$ there is a $d$ dimensional
quasiuniform mesh $\mcT_{h}(U_O)$ on $U_O$ with mesh parameter 
$h$, i.e. $\mcT_{h}(U_O)$ is a partition of $U_O$ into shape regular elements 
$T$ with diameter $h_T$ and $h_T \sim h$ uniformly. The active mesh is defined by 
\begin{equation}
\mcT_{h,O} = \{ T \in \mcT_{h}(U_O) : T \cap O \neq \emptyset \}
\end{equation}
which covers $O$ but is not required to perfectly match $O$. 

\item Let $V_{h,O}$ be a finite dimensional space consisting of continuous piecewise 
polynomial functions defined on $\mcT_{h,O}$ such that that there is an interpolation 
operator $\pi_{h,i}:H^1(O) \rightarrow V_{h,i}$, which satisfy the approximation property 
\begin{equation}\label{eq:interpol-first}
\| v - \pi_h v \|_{H^m(\mcT_{h,O})} \lesssim h^{s-m} \| v \|_{H^s(O)}\qquad m=0,1, 
\qquad m \leq s \leq p+1
\end{equation}
where $p+1$ is the approximation order of $V_{h,O}$. 

\item For $O \in \{\Omega_i\}_{i=1}^N$ we use the simplified notation 
\begin{equation}
 \mcT_{h,\Omega_i} = \mcT_{h,i}, \qquad V_{h,\Omega_i} = V_{h,i}
\end{equation}
we also define 
\begin{equation}
\mcT_{h,0}= \sqcup_{i=1}^{N_0} \mcT_{h,0,i}
\end{equation}
and for $O\in \{ \Omega_{0,i}\}_{i=1}^{N_0}$, 
\begin{equation}
\mcT_{h,\Omega_{0,i}} = \mcT_{h,0,i}, \qquad V_{h,\Omega_{0,i}} = V_{h,0,i}
\end{equation}

\item Define the finite element spaces
\begin{equation}\label{eq:Vh0-Vh1N}
V_{h,0} = \bigoplus_{i=1}^{N_0} V_{h,0,i},
\qquad 
V_{h,1,N} = \bigoplus_{i=1}^{N} V_{h,i}
\end{equation}
and
\begin{equation}\label{eq:Wh}
W_h  = V_{h,0} \oplus V_{h,1,N}
\end{equation}
\end{itemize}
See Figures \ref{fig:three-parts-matching-meshes} and \ref{fig:three-parts-nonmatching-meshes} in Section 4 for illustrations of the construction of the mesh.

\begin{rem} Here we are using restrictions of $d$ dimensional spaces to the $d-1$ dimensional 
skeleton in the spirit of CutFEM, see \cite{BHL15,ORG09}, where similar ideas are used to 
solve the Laplace--Beltrami problem on a codimension one surface embedded in $d$ dimensional 
mesh. We may also use a standard $d-1$ dimensional mesh on each of the components $\Omega_{0,i}$ 
of the skeleton. Note that these meshes does not need to match on the interfaces $\overline{\Omega}_{0,i} \cap \overline{\Omega}_{0,j}$ between two components of the skeleton.  Our analysis readily extends to 
this situation as well.
\end{rem}

\paragraph{Stabilization Forms.}
We here define a number of abstract properties, which we assume our stabilization forms satisfy. 
In Section~\ref{section:stabilization-forms} below we construct forms that fulfill these properties.
For $i=0,\dots,N$ there are symmetric bilinear forms $s_{h,i}$, on $V_{h,i} + H^{p+1}(\Omega_i)$, 
with associated semi norm $\| \cdot \|_{s_{h,i}}$ such that:
\begin{itemize}
\item For $v \in H^{p+1}(\Omega_i)$ it holds
\begin{equation}\label{eq:sh-consistency}
s_{h,i}(v,w) = 0 \qquad  \forall w \in V_{h,i} 
\end{equation}
and
\begin{equation}\label{eq:sh-approximation}
\| v - \pi_h v \|_{s_{h,i}} \lesssim h^p \| v \|_{H^s (\Omega_i)} 
\end{equation}
\item For $v \in V_{h,0,i}$, $i=1,\dots,N_0$ there are constants such that
\begin{equation}\label{eq:sh-inverseL2-Omega0}
h^{-1} \| v \|^2_{\mcT_{h,0,i}}
\sim
\| v \|^2_{\Omega_{0,i}} + \| v \|^2_{s_{h,0,i}} \qquad v \in V_{h,0,i} 
\end{equation}
\item
For $v \in V_{h,i}$, $i=1,\dots,N$, there is a constant such that
\begin{equation}\label{eq:sh-inverseL2}
\| v \|^2_{\mcT_{h,i}}
\lesssim 
\| v \|^2_{\Omega_i} + \| v \|^2_{s_{h,i}} \qquad v \in V_{h,i} 
\end{equation}
\item For $v \in V_{h,i}$, $i=1,\dots,N$, there is a constant such that
\begin{equation}\label{eq:sh-inversetrace}
h \| \nabla v \|^2_{\partial \Omega_i, a_i}
\lesssim 
\| \nabla v \|^2_{\Omega_i, a_i} + \| v \|^2_{s_{h,i}} \qquad v \in V_{h,i} 
\end{equation}
\end{itemize}

\paragraph{Derivation of the Method.}
Multiplying \eqref{eq:strong-a} by a test function $v$,
integrating by parts over $\Omega_i$, and using the convention 
$u_0 = 0$ on $\partial \Omega$ we obtain
\begin{align}\label{eq:der-a}
\sum_{i=1}^N (f_i,v_i)_\Omega 
=&{}\sum_{i=1}^N -(\nabla \cdot a_i \nabla  u_i, v_i)_{\Omega_i} 
\\
=&{}\sum_{i=1}^N (a_i \nabla u_i,\nabla v_i)_{\Omega_i} 
      - (n_i \cdot a_i \nabla u_i , v_i )_{\partial \Omega_i}
\\
=&{}\sum_{i=1}^N (a_i \nabla u_i,\nabla v_i)_{\Omega_i} 
      - (n_i \cdot a_i \nabla u_i , v_i )_{\partial \Omega_i}
      \\ \label{eq:der-b}
      & \qquad
      -(u_i-u_0,n\cdot a_i \nabla v_i)_{\partial \Omega_i} 
      +(\beta h^{-1} a_i (u_i - u_0), v_i)_{\partial \Omega_i} 
      \\ \label{eq:der-c}
=&{}\sum_{i=1}^N (a_i \nabla u_i,\nabla v_i)_{\Omega_i}  + s_{h,i}(u_i,v_i)
      - (n_i \cdot a_i \nabla u_i , v_i - v_0 )_{\partial \Omega_i}
      \\ 
      & \qquad
      -(u_i-u_0,n_i \cdot a_i \nabla v_i)_{\partial \Omega_i} 
      +(\beta h^{-1} a_i (u_i - u_0), (v_i - v_0) )_{\partial \Omega_i} 
\end{align}
where we added terms which are zero for the exact solution $u$ with the 
purpose of obtaining a stable symmetric formulation. More precisely, in 
(\ref{eq:der-b}) we used (\ref{eq:strong-c}) to add terms involving  
$[u]_i = u_i - u_0 = 0$ and in (\ref{eq:der-c}) 
we used (\ref{eq:strong-b}) and the identity 
\begin{equation}
0 = (\llbracket n \cdot a \nabla u \rrbracket, v_0)_{\Omega_0} 
= \sum_{i=1}^N (n_i \cdot  a_i\nabla u_i, v_0)_{\partial \Omega_i}
\end{equation}
since $v_0 = 0$ on $\partial \Omega$. Finally, we added the stabilization form 
$s_{h,i}$ using consistency (\ref{eq:sh-consistency}).

\paragraph{Definition of the Method.} Find $u_h \in W_h$ such that 
\begin{equation}\label{eq:fem}
A_h(u_h,v) = l_h(v)\qquad \forall v \in W_h
\end{equation}
where $W_h$ is defined in (\ref{eq:Wh}) and
\begin{align}
A_h(v,w) = &{} \sum_{i=1}^N (a_i \nabla v_i,\nabla w_i)_{\Omega_i} + s_{h,i}(v_i,w_i)
\\ & \qquad \label{eq:method-nitsche-a}
-(n_i \cdot a_i \nabla v_i, [w]_i)_{\partial \Omega_i}
-([v]_i,n_i \cdot a_i \nabla w_i)_{\partial \Omega_i}
\\ & \qquad \label{eq:method-nitsche-b}
+(\beta h^{-1}_i a_i [v]_i, [w]_i)_{\partial \Omega_i}
\\ &
+ s_{h,0}(v_0,v_0)
\\
l_h(v) &= \sum_{i=1}^N (f_i,v_i)_{\Omega_i}
\end{align}

\begin{rem} We consider the partition $\{\Omega_i\}_{i=1}^N$ of $\Omega$ 
as being fixed for simplicity only. In fact we may also allow the number of 
subdomains $N$ to increase during refinement.  We then obtain a version 
of the polygonal finite element method, where each polygonal finite element 
is equipped with a mesh and piecewise continuous polynomials of degree 
$p$, and not only polynomials of order $p$ which in general is the case in 
polygonal finite elements . In order to guarantee that coercivity holds under 
refinement of the partition we essentially need the  inverse inequality 
(\ref{eq:sh-inversetrace}) to hold with a uniform constant for all subdomains 
that occur during the refinement of the partition. Therefore  we  need to 
assume some additional shape regularity of the subdomains. For instance, 
we may assume that for each $\Omega_j$ there is partition $\mcS_j$ of 
$\Omega_i$ into simplicies $S$ and a constant such that for all 
$S \in \mcS_j$ with $|\partial S \cap \partial \Omega|>0$ it holds
\begin{equation}\label{eq:shape-reg}
h_{i} \lesssim \frac{|S|}{|\partial S \cap \partial \Omega_i|} 
\end{equation}  
where the mesh parameter $h_i \in (0,h_{0,i}]$. Then (\ref{eq:sh-inversetrace}) 
holds uniformly over all partitions $\{\Omega_j\}_{j=1}^N$ of $\Omega$, with 
corresponding meshes $\mcT_{h,i}$.

The corresponding condition for polygonal finite elements is 
\begin{equation}\label{eq:shape-reg-poly}
\text{diam}(\Omega_i) \lesssim \frac{|S|}{|\partial S \cap \partial \Omega_i|} 
\end{equation}  
see Assumption 30, Section 4.3 in \cite{CDGH17}, and since we may always 
assume that $h_{i} \lesssim  \text{diam}(\Omega_i)$ we conclude that 
(\ref{eq:shape-reg}) is weaker than (\ref{eq:shape-reg-poly}) and thus more 
complex subdomains may be employed when finer meshes are used.

\end{rem}

\section{Error Estimates}
The error estimates for our methods are obtained in a similar fashion
as those proven in \cite{Eg09} and we will focus here on how
robustness and optimal estimates are obtained in the framework of cut elements.
\subsection{Properties of $\boldsymbol{A}_{\boldsymbol{h}}$}

Define the energy norm 
\begin{equation}\label{eq:energy-norm}
\tn v \tn^2_h 
= 
\sum_{i=1}^N \| \nabla v_i \|^2_{\Omega_i,a_i}
+h \| \nabla v_i \|^2_{\partial \Omega_i,a_i}
+h^{-1} \| [ v ]_i\|^2_{\partial \Omega_i,a_i}
\end{equation}

\paragraph{Continuity.} It follows directly from the Cauchy-Schwarz inequality 
that $A_h$ is continuous
\begin{equation}\label{eq:continuity-Ah}
A_h(v,w) \lesssim \tn v \tn_h \tn w \tn_h\qquad  v,w \in W + W_h
\end{equation}

\paragraph{Coercivity.} For $\beta$ large enough $A_h$ is  coercive 
\begin{equation}\label{eq:coercivity-Ah}
\tn v \tn^2_h \lesssim A_h(v,v)\qquad v \in W_h
\end{equation}
\begin{proof}[Proof of (\ref{eq:coercivity-Ah}).]Using standard arguments and 
the  inverse estimate (\ref{eq:sh-inversetrace}) we obtain 
\begin{align}
A_h(v,v) 
= &{} \sum_{i=1}^N (a_i \nabla v_i, \nabla v_i)_{\Omega_i} + s_{h,i}(v,w)
- 2(n_i \cdot a_i \nabla v_i, [v_i])_{\partial \Omega_i} \\& \qquad
+\beta h^{-1} \| [v]_i\|^2_{\partial \Omega_i,a_i}
\\
\gtrsim &{}
\sum_{i=1}^N \|\nabla v_i \|^2_{\Omega_i,a_i} + \|v\|^2_{s_{h,i}}
- 2 h^{1/2} \| \nabla v_i\|_{\partial \Omega_i,a_i } 
h^{-1/2} \|[v_i]\|_{\partial \Omega_i,a_i} \\& \qquad
+\beta h^{-1} \| [v]_i\|^2_{\partial \Omega_i,a_i}
\\
\gtrsim  &{}
\sum_{i=1}^N \Big( \|\nabla v_i \|^2_{\Omega_i,a_i} + \|v\|^2_{s_{h,i}}
- \delta h \| \nabla v_i\|^2_{\partial \Omega_i,a_i} \Big) \\ & \qquad
+  
\bigl( \beta  - \delta^{-1} \bigr) h^{-1}  \| [v]_i\|^2_{\partial \Omega_i,a_i}
\\
\gtrsim  &{}
\sum_{i=1}^N \bigl(1 - \delta C_{1}\bigr) \Big(\|\nabla v_i \|^2_{\Omega_i,a_i} + \|v\|^2_{s_{h,i}}\Big)
+  
\bigl( \beta  - \delta^{-1} \bigr) h^{-1}  \| [v]_i\|^2_{\partial \Omega_i,a_i}
\end{align}
for $\delta>0$ sufficiently small. Here $C_1$ is the hidden constant in (\ref{eq:sh-inversetrace}).
\end{proof}
\paragraph{Poincar\'e Inequality.} Let $\phi$ 
be the solution to the dual problem 
\begin{equation}\label{eq:poincare-dual}
\text{$-\nabla \cdot a \nabla \phi = \psi$ in $\Omega$} ,
\qquad 
\text{$\phi = 0$ on $\partial \Omega$}
\end{equation} 
with $\psi \in L^2(\Omega)$, and assume that 
\begin{equation}\label{eq:poincare-dual-reg}
\sum_{i=1}^N \| \phi_i \|^2_{H^{3/2}(\Omega_i)} 
\lesssim 
\|\psi\|^2_\Omega
\end{equation}
Under the regularity assumption (\ref{eq:poincare-dual-reg})  the following Poincar\'e estimate holds 
\begin{align}\label{eq:poincare}
\Big( \min_{1\leq i \leq N} d_{\Omega_i} \Big) h^{-1} \| v_0 \|^2_{\mcT_{h,0}} +  \sum_{i=1}^N \| v_i \|^2_{\mcT_{h,i}}  \lesssim \tn v \tn^2_h  \qquad v \in W_h
\end{align}
where $d_{\Omega_i}$ is the diameter of $\Omega_i$. This estimate in particular shows that 
$\tn \cdot \tn_h$ is a norm on $W_h$.
\begin{proof}[Proof of  (\ref{eq:poincare}).] Consider first the  estimation of the bulk 
subdomain contributions we first note that using (\ref{eq:sh-inverseL2}) we have
\begin{equation}
\sum_{i=1}^N \| v \|^2_{\mcT_{h,i}} \lesssim \sum_{i=0}^N \| v \|^2_{\Omega_i} + \| v \|^2_{s_{h,i}}
\end{equation}
Next to estimate  the term $\sum_{i=1}^N \| v \|^2_{\Omega_i}$ we multiply the dual problem 
(\ref{eq:poincare-dual}) by $v \in W_h$ and using partial integration on each $\Omega_i$ 
we obtain 
\begin{align}\label{eq:poincare-proof-a}
 \sum_{i=1}^N (v_i,\psi_i)_{\Omega_i} 
 &=
 \sum_{i=1}^N ( a_i \nabla v_i, \nabla \phi_i )_{\Omega_i} - (v_i, n_i \cdot a_i \nabla \phi_i )_{\partial \Omega_i}
 \\ \label{eq:poincare-proof-b}
 &=
 \sum_{i=1}^N ( a_i \nabla v_i, \nabla \phi_i )_{\Omega_i} 
 - ([v]_i, n_i \cdot a_i \nabla \phi_i )_{\partial \Omega_i}
- (v_0, n_i \cdot a_i \nabla \phi_i )_{\partial \Omega_i}
 \\ \label{eq:poincare-proof-c}
 &=
 \sum_{i=1}^N ( a_i \nabla v_i, \nabla \phi_i )_{\Omega_i} 
 - ([v]_i, n_i \cdot a_i \nabla \phi_i )_{\partial \Omega_i} 
 - (v_0, \underbrace{\llbracket n \cdot a \nabla \phi \rrbracket}_{=0} )_{\Omega_0}
 \\ \label{eq:poincare-proof-d}
 &\lesssim
  \sum_{i=1}^N \| \nabla v_i\|_{\Omega_i,a_i} \|\nabla \phi_i \|_{\Omega_i,a_i} 
 + h^{-1/2} \|[v]_i\|_{\partial \Omega_i,a_i} h^{1/2} \| \nabla \phi_i \|_{\partial \Omega_i,a_i}
  \\ \label{eq:poincare-proof-e}
 &\lesssim
 \tn v \tn_h \underbrace{\Big( \sum_{i=1}^N \|\nabla \phi \|^2_{\Omega_i,a_i} 
 + h \| \phi \|^2_{H^{3/2}(\Omega_i),a_i} \Big)^{1/2}}_{ \lesssim  \|\psi \|_{\Omega} }
\\ \label{eq:poincare-proof-f}
&\lesssim \tn v \tn_h \| \psi \|_{\Omega}
\end{align}
where in (\ref{eq:poincare-proof-b}) we added and subtracted $v_0$ 
in the boundary terms; in (\ref{eq:poincare-proof-c}) we used the fact 
that $v_0= 0$ on $\partial \Omega$; in (\ref{eq:poincare-proof-d}) we 
used the Cauchy-Schwarz inequality and the fact that 
$\llbracket n\cdot a \nabla \phi \rrbracket =0$ 
on $\Omega_0$, in (\ref{eq:poincare-proof-e}) we used the definition 
(\ref{eq:energy-norm}) of the energy norm and a trace inequality for $\phi$;
and finally in (\ref{eq:poincare-proof-f}) we used the regularity assumption (\ref{eq:poincare-dual-reg}).
Thus we obtain 
\begin{equation}
\sum_{i=1}^N \| v_i \|^2_{\mcT_{h,i}} \lesssim \tn v \tn_h^2
\end{equation}
Next using the trace inequality 
\begin{equation}
\| v \|^2_{\partial \Omega_i} \lesssim d^{-1}_{\Omega_i} \| v_i \|^2_{\Omega_i} + d_{\Omega_i} 
\| \nabla v_i \|^2_{\Omega_i}
\end{equation}
we have that
\begin{align}
\| v_0\|^2_{\Omega_0} 
&\lesssim 
\sum_{i=1}^N \|[v]_i \|^2_{\partial \Omega_i}  + \| v_i \|^2_{\partial \Omega_i} 
\\
 &\lesssim 
\sum_{i=1}^N \|[v]_i \|^2 + d^{-1}_{\Omega_i} \| v_i \|^2_{\Omega_i} + d_{\Omega_i} 
\| \nabla v_i \|^2_{\Omega_i}
\\
&\lesssim \Big( \min_{1\leq i \leq N} d_{\Omega_i} \Big)^{-1} \tn v \tn_h^2 
\end{align}
and thus by (\ref{eq:sh-inverseL2-Omega0}) we obtain
\begin{align}
h^{-1} \| v_0 \|^2_{\mcT_{h,0}} 
\lesssim 
\| v_0 \|^2_{\Omega_0} + \| v_0 \|^2_{s_{h,0}} 
\lesssim 
 \Big( 1 + \Big( \min_{1\leq i \leq N} d_{\Omega_i} \Big)^{-1} \Big ) \tn v \tn_h^2
\end{align}
which concludes the proof since  $\bigl( \min_{1\leq i \leq N} d_{\Omega_i} \bigr)^{-1} \lesssim 1$.
\end{proof}

Note in particular that we have the estimate 
\begin{equation}\label{eq:poincare-omega0}
\Big( \min_{1\leq i \leq N} d_{\Omega_i} \Big) h^{-1} \| v_0 \|^2_{\mcT_{h,0}}
\lesssim \tn v_0 + w \tn_h
\end{equation}
for all $w \in V_{h,1,N} = \oplus_{i=1}^N V_{h,i}$.

\subsection{Interpolation Operator}
\paragraph{Definition of the Interpolation Operator.}
\begin{itemize}
\item In order to handle both the $d$ dimensional components and the $d-1$ dimensional 
components at the same time we let $O \in \{ \Omega_{0,i}\}_{i=1}^{N_0} \cup \{\Omega_i\}_{i=1}^N$, 
i.e., $O$ is either one of the $d$ dimensional bulk domains $\Omega_i$ or one of 
the $d-1$ dimensional components of the skeleton $\Omega_0$.

\item There is an extension operator 
$E_O:H^s(O) \rightarrow U_{\delta}(O)$ such that 
\begin{equation}\label{ext:stab}
\| E v \|_{H^s(U_\delta(O))} \lesssim \delta^{1/2} \| v \|_{H^s(O)}
\end{equation}
where $U_\delta(O) = \cup_{x\in O} B_\delta(x)$, with $B_\delta(x)$
denoting a $d$ dimensional ball with radius $\delta>0$ centered in $x$. Note 
that this means that $U_\delta(O)$ is always a $d$ dimensional tubular 
neighborhood of $O$.

\item Let $\mcT_h(O)$ be the mesh associated with $O$ and $V_h(O)$ 
the correponding finite element space.  Let $\pi_{h,O,Cl} : L^2(\mcT_h(O)) \rightarrow V_{h}(O)$ 
be the Clement interpolant and define $\pi_{h,O} : H^s(O) \rightarrow V_{h}(O)$ 
such that 
\begin{equation}
\pi_{h,O} v = \pi_{h,O,Cl} E_O v
\end{equation} 
Finally, we let $\pi_h: \left( \oplus_{i=1}^{N_0} H^s(\Omega_i)\right) 
\oplus \left(  \oplus_{i=1}^N H^s(\Omega_i) \right) \rightarrow V_h$ be 
defined componentwise by $(\pi_h v)_O = \pi_{h,O} v_O$. We have the 
interpolation error estimates 
\begin{equation}\label{eq:interpol-basic}
\| v - \pi_h v \|_{H^m(O)} \lesssim h^{s-m} \| v \|_{H^{s}(O)}\qquad m=0,1
\end{equation}
for $0\leq m \leq s \leq p+1$. 
For the $d$ dimensional domains $O=\Omega_i$, $i=1,\dots,N$, we derive 
(\ref{eq:interpol-basic}) using a standard interpolation error bound for the 
Clement interpolant followed by the stability of the extension operator
\begin{align}\label{eq:interpol-basic-prf-a}  
\| v - \pi_h E v \|_{H^m(\Omega_i)} 
&\leq 
\| v - \pi_h E v \|_{H^m(\mcT_{h,i})} 
\\ \label{eq:interpol-basic-prf-b} 
&\qquad \leq 
 h^{p+1-m} \| E v \|_{H^{p+1}(\mcT_{h,i})}
\leq 
 h^{p+1-m} \| v \|_{H^{p+1}(\Omega_i)}
\end{align}
For the $d-1$ dimensional skeleton subdomains $\Omega_{0,i}$, $i=1,\dots, N_0$, we instead 
employ the trace inequality 
\begin{equation}\label{eq:trace}
\| v \|^2_{\Omega_{0,i}} \lesssim h^{-1} \| v \|^2_{\mcT_{h,0,i}} + h \| \nabla v \|^2_{\mcT_{h,0,i}}  
\end{equation}
and then (\ref{eq:interpol-basic}) can be derived in the same way as in 
(\ref{eq:interpol-basic-prf-a})-(\ref{eq:interpol-basic-prf-b}).
\end{itemize}

\paragraph{Interpolation Error Estimate.} The following estimate holds
\begin{equation}\label{eq:interpol-energy}
\tn v - \pi_h v \tn^2_h 
\lesssim  h^{2p}  \| v_0 \|^2_{H^{p+1/2}(\Omega_0)} 
+ \sum_{i=0}^N   h^{2p} \| v_i \|^2_{H^{p+1}(\Omega_i)} 
\end{equation}

\begin{rem}
Note that here we can easily localize $h$ and $p$ to the subdomains and skeleton subdomains. But for clarity we keep a global mesh parameter $h$ and order of polynomials $p$.
\end{rem}

\begin{proof}[Proof of (\ref{eq:interpol-energy}).] Let $\eta = v - \pi_h v$ be the interpolation 
error and using the definition (\ref{eq:energy-norm}) of the energy norm we have
\begin{align}
\tn \eta \tn^2_h 
&= \sum_{i=1}^N \| \nabla \eta \|^2_{\Omega_i} 
+ h \| \nabla \eta \|^2_{\partial \Omega_i} + h^{-1} \| [\eta]_i \|^2_{\partial \Omega_i} 
\\
&\leq\sum_{i=1}^N \| \nabla \eta \|^2_{\Omega_i} + h \| \nabla \eta \|^2_{\partial \Omega_i} 
+ h^{-1} \| \eta_i \|^2_{\partial \Omega_i}  + h^{-1} \| \eta_0 \|^2_{\partial \Omega_i} 
\\
&\lesssim \sum_{i=1}^N \| \nabla \eta \|^2_{\Omega_i} + h \| \nabla \eta \|^2_{\partial \Omega_i} 
+ h^{-1} \| \eta_i \|^2_{\partial \Omega_i}  + \sum_{i=1}^{N_0} h^{-1} \| \eta_0 \|^2_{\Omega_{i,0}} 
\\
&\lesssim \sum_{i=1}^N \| \nabla \eta \|^2_{\Omega_i} 
+ ( \| \nabla \eta \|^2_{\mcT_{h,i}(\partial \Omega_i)} 
+ h^2 \| \nabla^2 \eta \|^2_{\mcT_{h,i}(\partial \Omega_i)} )
\\ 
&\qquad \qquad + ( h^{-2} \|  \eta \|^2_{\mcT_{h,i}(\partial \Omega_i)}  
+ \| \nabla \eta \|^2_{\mcT_{h,i}(\partial \Omega_i)} ) 
+ \sum_{i=1}^{N_0} 2h^{-1} \| \eta_0 \|^2_{\Omega_{i,0}} 
\\
&\lesssim 
\sum_{i=1}^N h^{2p}  \| v  \|^2_{H^{p+1}(\Omega_i)} 
+ \sum_{i=1}^{N_0} h^{2p}  \| v_{0,i} \|^2_{H^{p+1/2}(\Omega_{0,i})}
\end{align}
Here we used the trace inequality 
\begin{equation}\label{eq:trace-partialomega}
\| v \|^2_{\partial \Omega_{i}} \lesssim h^{-1} \| v \|^2_{\mcT_{h,i}(\partial \Omega_i)} 
+ h \| \nabla v \|^2_{\mcT_{h,i}(\partial \Omega_i)}  
\end{equation}
where $\mcT_{h,i}(\partial \Omega_i)$ is the set of all elements in $\mcT_{h,i}$ which 
intersect $\partial \Omega_i$, and the interpolation estimate (\ref{eq:interpol-basic}).
\end{proof}

\subsection{Stabilization Forms} \label{section:stabilization-forms}
Define the following stabilization forms:
\begin{itemize}
\item For $i=1,\dots,N$ define the stabilization forms 
\begin{align}
s_{h,i}(v,w) = \sum_{l=1}^p c_{d,l} \, h^{2 l - 1} ([D_n^l v],[D_n^l w])_{\mcF_{h,i}}
\end{align}
where $c_{d,l}>0$ is a parameter and $\mcF_{h,i}$ is the set of interior faces in $\mcT_h$ which belongs 
to an element that intersect the boundary $\partial \Omega_i$ and 
$D^l_\xi v$ denotes the contraction of $\xi^l$ and $\nabla^l v$.

\item For $i=0$ recall that $\Omega_0= \cup_{j=0}^{N_0} \Omega_{0,j}$ 
and define 
\begin{equation} \label{eq:stab-skeleton}
s_{h,0} = \sum_{j=1}^{N_0} s_{h,0,j}
\end{equation}
where 
\begin{equation} \label{eq:stab-skeleton-sub}
s_{h,0,j}(v,w) =  \sum_{l=1}^p c_{d-1,l} \, h^{2l}\left( (D^l_\nu v, D^l_\nu w )_{\Omega_{0,j}} 
+  
([D_n^l v],[D_n^l w])_{\mcF_{h,0,j}} \right)
\end{equation}
where $c_{d-1,l}>0$ is a parameter,
$\nu$ is the normal to the smooth $d-1$ manifold $\Omega_{0,j}$
and $\mcF_{h,0,j}$ denotes the set of interior faces in $\mcT_{h,0,j}$.
\end{itemize}
Then the required properties of the stabilization forms 
(\ref{eq:sh-approximation})--(\ref{eq:sh-inversetrace}) hold, see 
\cite{MassingLarsonLoggEtAl2013a,HaLaLa18,LaZa2017}.

\subsection{Error Estimates}

\paragraph{Energy Norm Error Estimate.} The following estimate holds
\begin{equation} 
\tn u - u_h \tn^2_h \lesssim  h^{2p}  \| u_0 \|^2_{H^{p+1/2}(\Omega_0)}  
+  \sum_{i=0}^N   h^{2p} \| u_i \|^2_{H^{p+1}(\Omega_i)} 
\end{equation}
This result follows immediately from the coercivity, Galerkin orthogonality, 
continuity, and interpolation error estimate. For the readers
convenience we detail the proof.
\begin{proof} Let $e_h = \pi_h u - u_h$ and $\rho = u - \pi_h u$. Then we have
\begin{align}
\tn e_h \tn_h^2 &\lesssim A_h(e_h,e_h)
\\
&=A_h(\rho, e_h)
\\
&\leq \tn \rho \tn_h \tn e_h \tn_h
\\
&\lesssim \left( h^{2p}  \| u_0 \|^2_{H^{p+1/2}(\Omega_0)}  
+  \sum_{i=0}^N   h^{2p} \| u_i \|^2_{H^{p+1}(\Omega_i)} \right)^{1/2}  \tn e_h \tn_h
\end{align}
\end{proof}

\paragraph{Subdomain $\boldsymbol{L}^{\boldsymbol{2}}$  Error Estimates.} Assuming that the solution to the 
dual problem
\begin{equation}
\text{$-\nabla \cdot a \nabla \phi = \psi$ in $\Omega$},
\qquad
\text{$\phi = 0$ on $\partial \Omega$}
\end{equation}
has the regularity 
\begin{equation}
\sum_{i=1}^N \| \phi_i \|^2_{H^s(\Omega_i)} \lesssim \|\psi \|^2_\Omega, \qquad s \in (3/2,1]
\end{equation}
we have the following error estimate
\begin{equation}
\sum_{i=1}^N \| u_i - u_{h,i} \|_{\Omega_i}^2  \lesssim  h^{2p+2s}  \| u_0 \|^2_{H^{p+1/2}(\Omega_0)}  
+  \sum_{i=0}^N   h^{2p+2s} \| u_i \|^2_{H^{p+1}(\Omega_i)} 
\end{equation}

\paragraph{Skeleton $\boldsymbol{L}^{\boldsymbol{2}}$  Error Estimates.} In order to 
estimate the $L^2$ norm of the skeleton error we  instead consider the dual problem 
\begin{alignat}{3}\label{eq:strong-dual-a}
-\nabla \cdot a_i \nabla \phi_i &= 0 &\qquad &\text{in $\Omega_i$}
\\ \label{eq:strong-dual-b}
\llbracket n\cdot a \nabla \phi \rrbracket &= \psi & \qquad &\text{on $\Omega_0$}
\\ \label{eq:strong-dual-c}
[\phi]_i &=0  & \qquad &\text{on $\partial \Omega_i \cap \Omega_0$}
\\ \label{eq:strong-dual-d}
\phi_i&=0 & \qquad  &\text{on $\partial \Omega_i \cap \partial \Omega$}
\end{alignat}
since $\psi$ is in $L^2$ on the skeleton the maximal regularity is 
\begin{equation}
\sum_{i=1}^N \| \phi_i \|^2_{H^{3/2}(\Omega_i)} 
\lesssim \|\psi \|^2_{\Omega_0}
\end{equation}
and thus we arrive at the estimate
\begin{equation}
 \| u_0 - u_{h,0} \|^2  \lesssim  h^{2p+1}  \| u_0 \|^2_{H^{p+1/2}(\Omega_0)}  
+  \sum_{i=0}^N   h^{2p+1} \| u_i \|^2_{H^{p+1}(\Omega_i)} 
\end{equation}
Observe that this shows that the error in the discrete $H^{1/2}$ norm,
$h^{-\frac12} \| u_0 - u_{h,0} \|$ has similar convergence order as
the energy norm error, which is optimal.

\section{The Schur Complement}

In this section we show how the bulk fields can be eliminated using the Schur complement and we derive a bound for the condition number of stiffness matrix associated with the Schur complement.

\subsection{Definitions}
\begin{itemize}
\item Define the operator $T_{h} : V_{h,0} \rightarrow V_{h,1,N} = \bigoplus_{i=1}^N V_{h,i}$ 
such that 
\begin{equation}\label{eq:Thdef}
A_h(v_0 + T_h v_0 ,  w ) = 0 \qquad \forall w \in V_{h,1,N}
\end{equation}

\item Define the Schur complement form on  $V_{h,0}$ by 
\begin{equation}
S_h(v, w )_{\Omega_0} = A_h(v + T_h v, w + T_h w )\qquad v,w \in V_{h,0}  
\end{equation}

\item Define the energy norm on $V_{h,0}$ associated with the Schur complement 
 by 
\begin{equation}
\tn v_0 \tn_{S_h} = \tn v_0 + T_h v_0 \tn_h 
\end{equation}

\item It follows directly from the definition of the Schur complement form that $S_h$ 
is coercive and continuous on $V_{h,0}$.

\end{itemize}

\subsection{Solution Using The Schur Complement}
Note that we have the $A_h$-orthogonal splitting
\begin{equation}\label{eq:orth-split}
W_h = (I + T_h) V_{h,0} \perp ( \{0\} \oplus V_{h,1,N} )
\end{equation}
and thus $u_h = (I+T_h) u_{h,0} + (0 \oplus u_{h,1,N})$ where 
$u_{h,0} \in V_{h,0}$ is the solution to
\begin{equation}\label{eq:shur-prob}
S_h(u_{h,0}, v_0 ) = l_h(  ( I+ T_h ) v_0 ) \qquad \forall v_0 \in V_{h,0}
\end{equation} 
and $u_{h,1,N}$ is the solution to 
\begin{equation}\label{eq:subdomain-prob}
A_h(0 \oplus u_{h,1,N}, 0 \oplus v) = l_h(v) \qquad \forall v \in V_{h,1,N}
\end{equation}
where we used the notation $0 \oplus u_{h,1,N}$ to indicate that the 
component in $V_{h,0}$ is zero. We note that (\ref{eq:subdomain-prob}) 
decouples and can be solved subdomain wise. For each subdomain 
we get a Nitsche type formulation with homogeneous Dirichlet data 
and right hand side $f_i = f |_{\Omega_i}.$

\subsection{Some Basic Estimates}
We collect the basic bounds needed to show an estimate of the condition number of 
the stiffness matrix associated with the Schur complement.

\paragraph{Norm Equivalence.}
There are constants such that
\begin{equation}\label{eq:schur-eqv}
\inf_{w \in V_{h,1,N}} \tn v_0 + w \tn_h \lesssim 
\tn v_0 \tn_{S_h} \lesssim \inf_{w \in V_{h,1,N}} \tn v_0 + w \tn_h
\end{equation}
\begin{proof}[Proof of (\ref{eq:schur-eqv}).] {\bf 1.} By definition 
\begin{equation}
\inf_{w \in V_{h,1,N}} \tn v_0 + w \tn_h \lesssim \tn v_0 + T_h v_0 \tn_h = \tn v_0 \tn_{S_h}
\end{equation}
\paragraph{2.} We have
\begin{align}
\tn v_0 \tn^2_{S_h} &=  \tn v_0 + T_h v_0 \tn^2_h 
\\
&\lesssim A_h( v_0 + T_h v_0, v_0 + T_h v_0) 
\\
&=  A_h( v_0 + T_h v_0, v_0 + w)
\\
&\lesssim 
\tn v_0 + T_h v_0 \tn_h \tn v_0 + w \tn_h
\end{align}
where we used the coercivity (\ref{eq:coercivity-Ah}) of $A_h$, the orthogonality 
(\ref{eq:Thdef}) of $T_h$, and the continuity (\ref{eq:continuity-Ah}) of $A_h$.  
Thus we conclude that 
\begin{equation}
\tn v_0 \tn_h \lesssim \tn v_0 + w \tn_h 
\end{equation}
for all $w \in V_{h,1,N}$, and therefore 
\begin{equation}
\tn v_0 + T_h v_0 \tn_h  \lesssim   \inf_{w \in V_{h,1,N}} \tn v_0 + w \tn_h 
\end{equation}
\end{proof}

\paragraph{Skeleton Poincar\'e Estimate.}
There is a constant such that for all $v \in V_{h,0}$,  
\begin{equation}\label{eq:Poincare-special}
 \Big( \min_{1\leq i \leq N} d_{\Omega_i}\Big)  h^{-1} \| v \|^2_{\mcT_{h,0}} \lesssim \tn v \tn^2_{S_h}
\end{equation} 
\begin{proof}[Proof of (\ref{eq:Poincare-special}).] Using the Poincar\'e inequality 
(\ref{eq:poincare-omega0}) we have for any $w \in V_{h,1,N}$, 
\begin{equation}\label{eq:Poincare-special-a}
 \Big( \min_{1\leq i \leq N} d_{\Omega_i}\Big) h^{-1} \| v \|^2_{\mcT_{h,0}} 
 \lesssim \tn v + w \tn^2_h
\end{equation} 
and as a consequence
\begin{equation}\label{eq:Poincare-special-b}
 \Big( \min_{1\leq i \leq N} d_{\Omega_i}\Big) h^{-1} \| v \|^2_{\mcT_{h,0}} 
 \lesssim \tn v + T_h v  \tn^2_h 
 = \tn v \tn^2_{S_h}
\end{equation} 
\end{proof}

\paragraph{Skeleton Inverse Estimate.}
There is a constant such that for all $v \in V_{h,0}$, 
\begin{equation}\label{eq:Inverse-schur}
\tn v \tn^2_{S_h} \lesssim h^{-1} ( \| v \|^2_{\Omega_0} + \| v \|^2_{s_{h,0}} )
\end{equation}
\paragraph{Proof of (\ref{eq:Inverse-schur}).} We have 
\begin{align}
\tn v_0 \tn^2_{S_h} &\lesssim \inf_{w \in V_{h,1,N}} \tn v_0 + w \tn^2_h 
\\
&\lesssim \tn v_0 + 0 \tn^2_h
\\
&=\sum_{i=1}^N h^{-1} \| v_0 \|^2_{\partial \Omega_i}  + h^{-1} \| v \|^2_{s_{h,0}}
\\
&\lesssim h^{-1} ( \| v \|^2_{\Omega_0} + \| v \|^2_{s_{h,0}} )
\end{align}
where we choose $w = 0$ on $V_{h,1,N}$ and then used the definition of the energy 
norm.
\hfill $\blacksquare$

\subsection{Condition Number Estimate for the Schur Complement}
\paragraph{Definitions and Basic Results.}
\begin{itemize}
\item Let $\{\varphi_i\}_{i=1}^D$ be the basis in $V_{h,0} $ and denote the expansion 
by 
\begin{equation}
v  = \sum_{i=1}^D \widehat{v}_i \varphi_i
\end{equation} 
Then we have the equivalence 
\begin{equation}\label{eq:RnL2eqv}
\| v \|^2_{\mcT_{h,0}} \sim h^d \| \widehat{v}\|^2_{\IR^D}
\end{equation}

\item The stiffness matrix associated with the Schur complement is defined by 
\begin{equation}
(\widehat{S} \widehat{v},\widehat{w})_{\IR^D} = S_h( v, w ) 
\end{equation}
Then $\widehat{S}$ is SPD and thus the spectrum is real and positive, and we have 
the Rayleigh quotient characterization of the largest and smallest eigenvalues
\begin{align}\label{eq:Rayleigh}
\lambda_{\max} 
= 
\max_{\widehat{v} \in \IR^D \setminus \{0\}} 
\frac{(\widehat{S} \widehat{v},\widehat{v})_{\IR^D}}{\| \widehat{v}\|^2_{\IR^D}},
\qquad
\lambda_{\min} 
= 
\min_{\widehat{v} \in \IR^D \setminus \{0\}} 
\frac{(\widehat{S} \widehat{v},\widehat{v})_{\IR^D}}{\| \widehat{v}\|^2_{\IR^D}}
\end{align}

\item The condition number is defined by
\begin{equation}
\kappa = \frac{\lambda_{\max}}{\lambda_{\min}}
\end{equation}
where $\lambda_{\max}$ and $\lambda_{\min}$ denote the maximum and minimum 
eigenvalues of $\widehat{S}$. In view of (\ref{eq:Rayleigh}) we have the identity 
\begin{equation}
\kappa = \max_{\widehat{v} \in \IR^D \setminus \{0\}} 
\frac{(\widehat{S} \widehat{v},\widehat{v})_{\IR^D}}{\| \widehat{v}\|^2_{\IR^D}} 
\max_{\widehat{v} \in \IR^D \setminus \{0\}} 
\frac{\| \widehat{v}\|^2_{\IR^D}}{(\widehat{S} \widehat{v},\widehat{v})_{\IR^D}}
\end{equation}
\end{itemize}

\paragraph{Condition Number Estimate.} The condition number satisfies the estimate
\begin{equation}
\kappa(\widehat{S}) \lesssim h^{-1}  \Big( \min_{1\leq i \leq N} d_{\Omega_i}\Big)^{-1}
\end{equation}
\begin{proof} The bound follows directly from the following two estimates:

\noindent {\bf 1.} We have  
\begin{align}
(\widehat{S}\widehat{v},\widehat{v})_{\IR^D} 
&= A_h(v + T_h v, v + T_h v ) 
\\
&\lesssim \tn v + T_h v \tn^2_h 
\\
&\lesssim h^{-1} (\|v \|^2_{L^2(\Omega_0)} + \| v \|^2_{s_{h,0}} )
\\
&\lesssim h^{-2} \|v \|^2_{\mcT_{h,0} }
\\
&\lesssim h^{-2} h^{d} \| \widehat{v} \|^2_{\IR^D}
\end{align}
where we used the continuity (\ref{eq:continuity-Ah})  
of $A_h$, the inverse inequality (\ref{eq:Inverse-schur}), 
and the equivalences (\ref{eq:RnL2eqv}) and (\ref{eq:sh-inverseL2-Omega0}). Thus we conclude that 
\begin{equation}
\lambda_{\max}  \lesssim h^{d-2}
\end{equation}
{\bf 2.} We have
\begin{align}
(\widehat{S}\widehat{v},\widehat{v})_{\IR^D} 
&= A_h(v + T_h v, v + T_h v ) 
\\
&\gtrsim \tn v + T_h v \tn^2_h 
\\
&\gtrsim \Big( \min_{1\leq i \leq N} d_{\Omega_i}\Big) h^{-1}  \|v \|^2_{L^2(\mcT_{h,0})}
\\
&\gtrsim   \Big( \min_{1\leq i \leq N} d_{\Omega_i}\Big) h^{-1} h^{d} \| \widehat{v} \|^2_{\IR^D}
\end{align}
where we used the coercivity (\ref{eq:coercivity-Ah}) of $A_h$, 
the Poincar\'e estimate (\ref{eq:Poincare-special}), and 
the equivalence (\ref{eq:RnL2eqv}). We conclude that 
\begin{equation}
\lambda^{-1}_{\min}  \lesssim   h^{-(d-1)}  \Big( \min_{1\leq i \leq N} d_{\Omega_i}\Big)^{-1}
\end{equation}
\end{proof}


\section{Numerical Results}

For assessment and illustration we implemented a 2D version of the method. We first give some details on implementation choices and then we present some illustrating examples and convergence results.

\subsection{Implementation} \label{section:implemention}

\paragraph{Approximation Spaces.}
On each subdomain and skeleton subdomain we define an approximation space by a mesh equipped with some finite elements. In all examples below we use standard Lagrange elements of degree $p$ which on quadrilaterals are tensor product polynomial elements $Q_p$ and on triangles are full polynomial elements $P_p$.
While the mesh on each subdomain and skeleton subdomain could be constructed completely independent of each other we here focus on two cases:
\begin{itemize}
\item \emph{Global Background Grid.\ } All meshes, i.e. subdomain meshes $\{\mcT_{h,i}\}_{i=1}^N$ and skeleton subdomain meshes $\{\mcT_{h,0,k}\}_{k=1}^{N_0}$, are extracted from the same background grid $\mcT_{h,\Omega}$.

\item \emph{Single Element Interfaces.\ } Subdomain meshes $\{\mcT_{h,i}\}_{i=1}^N$ may be arbitrarily constructed while each skeleton subdomain $\Omega_{0,k}$, $1 \leq k \leq N_0$, is equipped with a mesh $\mcT_{h,0,k}$ consisting of a single, typically higher order, element in $\mathbb{R}^d$.
\end{itemize}
The benefit of either choice is that the implementation of the Nitsche terms \eqref{eq:method-nitsche-a}--\eqref{eq:method-nitsche-b} becomes particularly straightforward as we for $T \in \mcT_{h,i}$, $\partial\Omega_i \cap T \neq \emptyset$, trivially know the element of the neighboring skeleton subdomain mesh and that $\partial\Omega_i \cap T$ is completely contained within that element. More general cases require the construction of the union of subdomain meshes and skeleton subdomain meshes to correctly evaluate the skeleton integrals. For curved skeletons, parametrically mapped subdomains or simply the case $d=3$ this construction is typically challenging to perform in a robust and efficient manner.

\paragraph{Parameter Choices.}
For the Nitsche penalty parameter we choose $\beta = 10 \cdot p^2$
and for the stabilization parameters we, as in the numerical examples in \cite{LaZa2017}, choose
$c_{d,l} = c_{d-1,l} = \frac{10^{-3}}{l!}$. We note no particular sensitivity in the choice of values for these parameters.

\subsection{Numerical Examples}

For our numerical examples we consider different partitions of the unit square $[0,1]^2$. We let the right hand side be given by $f=1$ and we vary the material coefficient $a$ in each subdomain.

\paragraph{Example 1: Three Subdomains.}

\begin{figure}\centering
\begin{subfigure}[t]{0.3\linewidth}\centering
\includegraphics[width=0.9\linewidth]{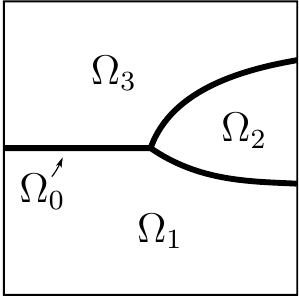}
\end{subfigure}
\caption{Illustrations of the model problem geometry. The unit square $[0,1]^2$ is divided into three subdomains according to the figure with material coefficients $a_1=1$, $a_2=2$ and $a_3=3$.}
\label{fig:model-problem}
\end{figure}

\begin{figure}\centering
\begin{subfigure}[t]{0.35\linewidth}\centering
\includegraphics[width=0.9\linewidth]{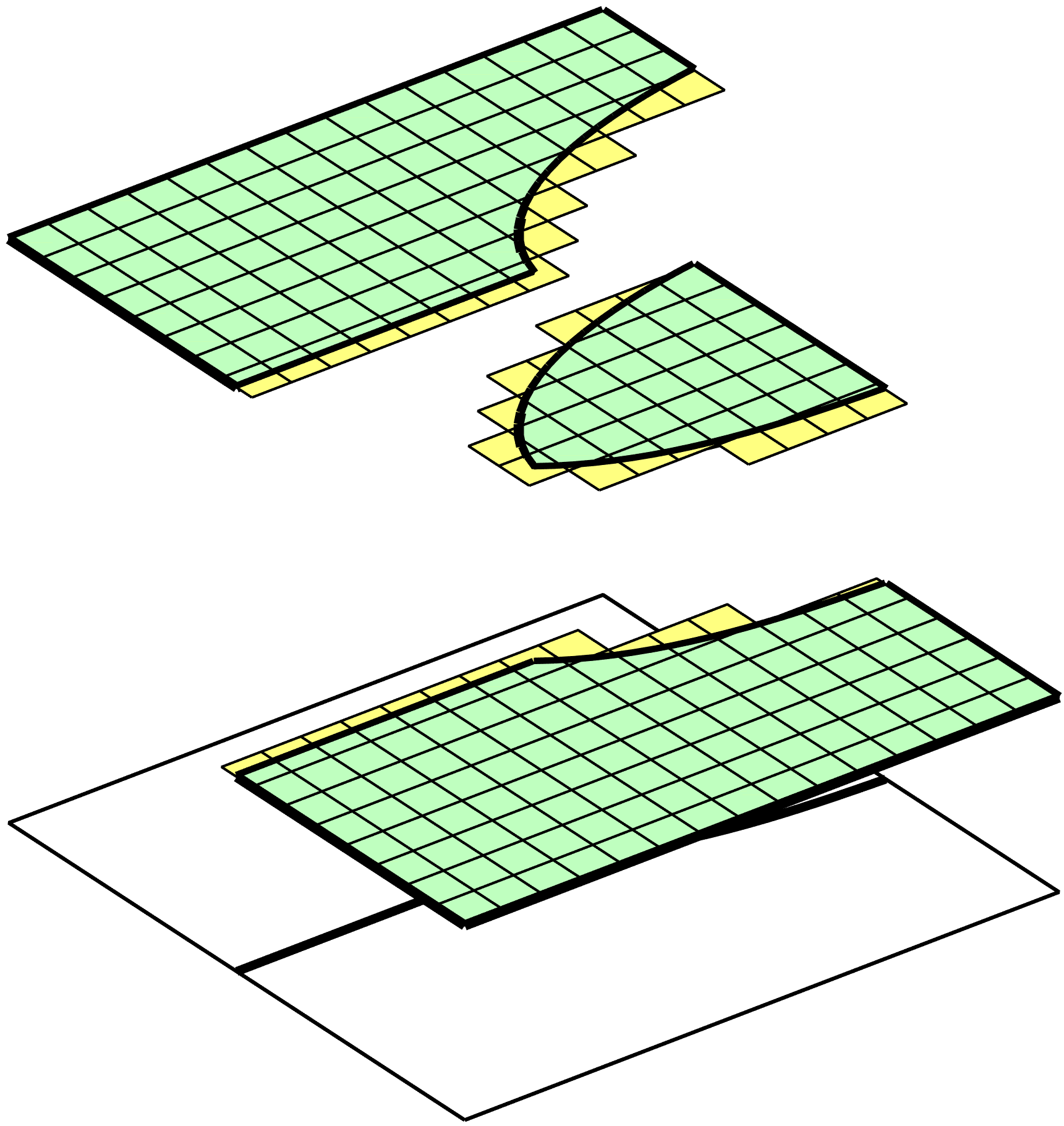}
\subcaption{Subdomain meshes}
\end{subfigure}
\begin{subfigure}[t]{0.35\linewidth}\centering
\includegraphics[width=0.9\linewidth]{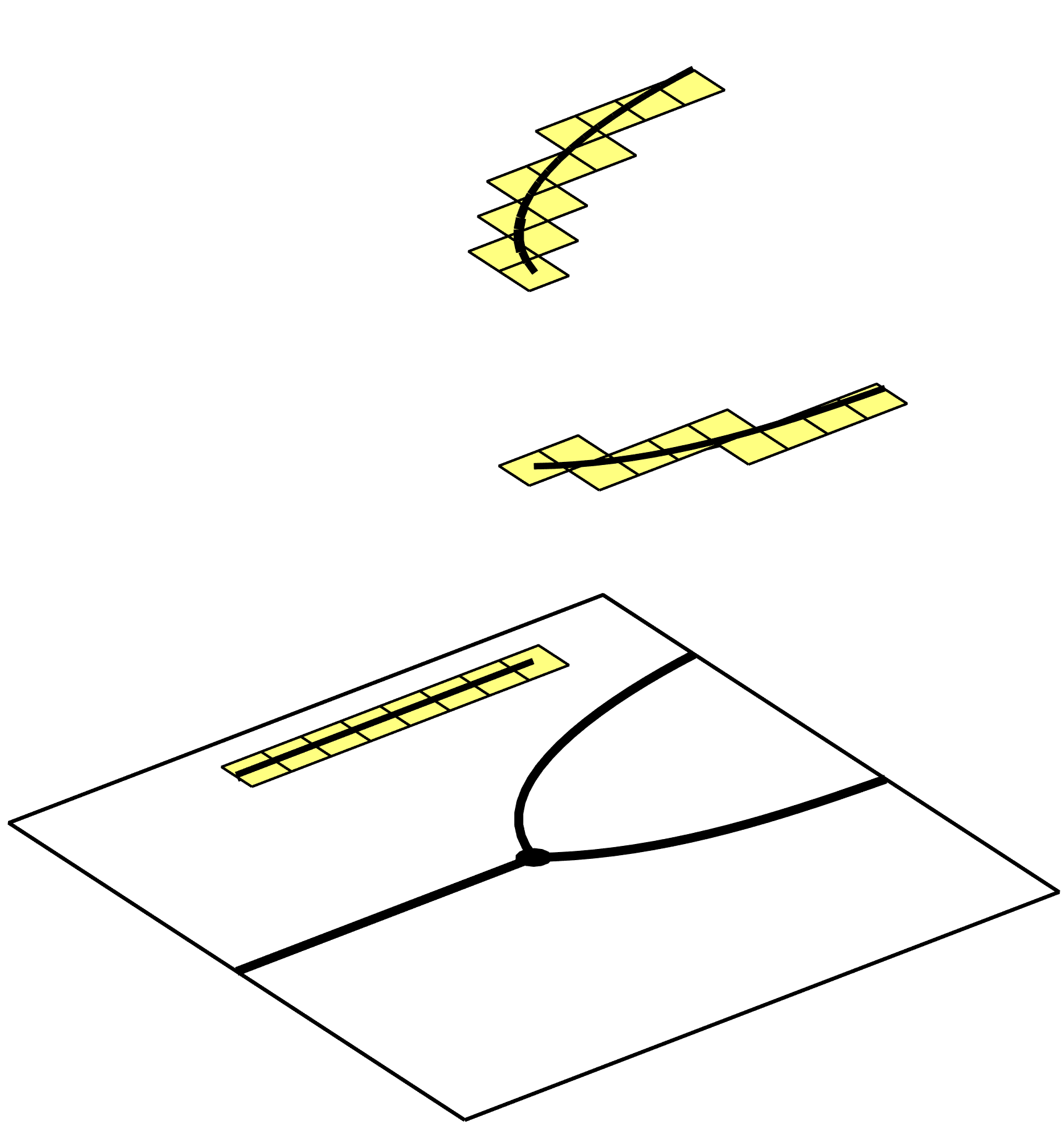}
\subcaption{Skeleton subdomain meshes}
\end{subfigure}
\caption{Meshes in the three subdomains example extracted from a global background grid. (a) Subdomain solutions are approximated using $Q_2$ elements. (b) Skeleton subdomain solutions are approximated in an embedding space of $Q_2$ elements.}
\label{fig:three-parts-matching-meshes}
\end{figure}

\begin{figure}\centering
\begin{subfigure}[t]{0.40\linewidth}\centering
\includegraphics[width=0.9\linewidth]{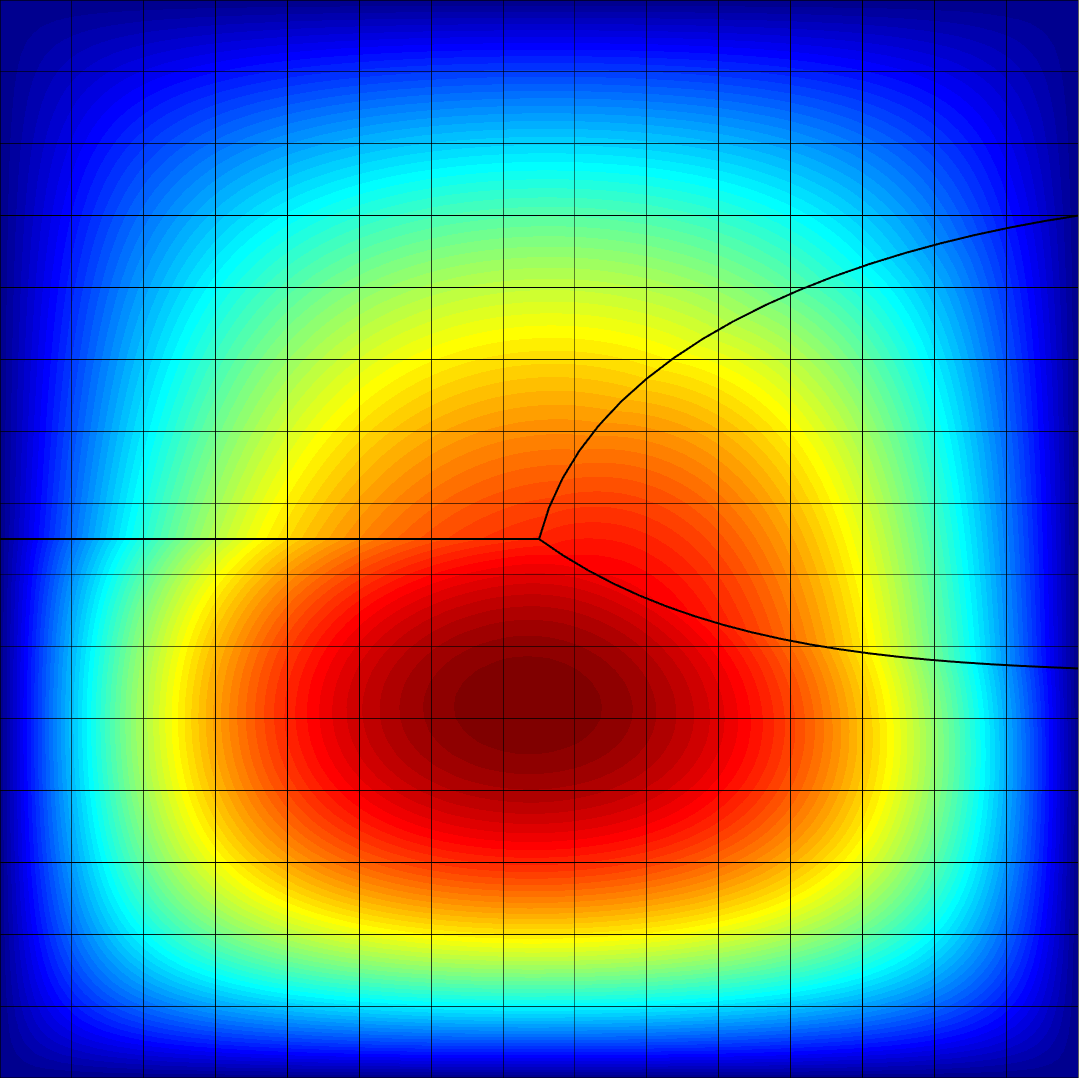}
\subcaption{Solution}
\end{subfigure}
\begin{subfigure}[t]{0.40\linewidth}\centering
\includegraphics[width=0.9\linewidth]{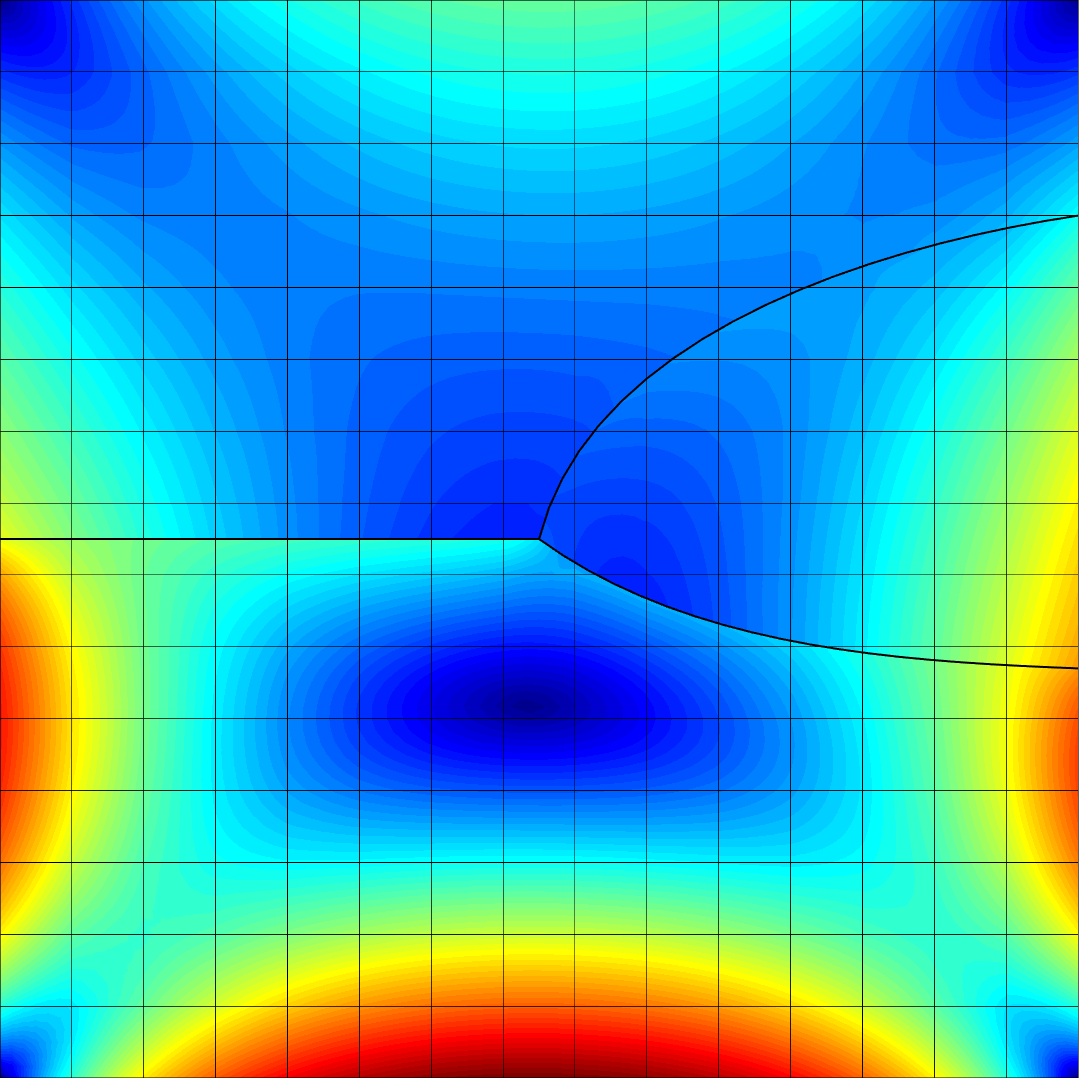}
\subcaption{Gradient magnitude}
\end{subfigure}
\caption{Three subdomains with different material coefficients and $Q_2$ meshes extracted from a global background grid.}
\label{fig:three-parts-matching}
\end{figure}

\begin{figure}\centering
\begin{subfigure}[t]{0.35\linewidth}\centering
\includegraphics[width=0.9\linewidth]{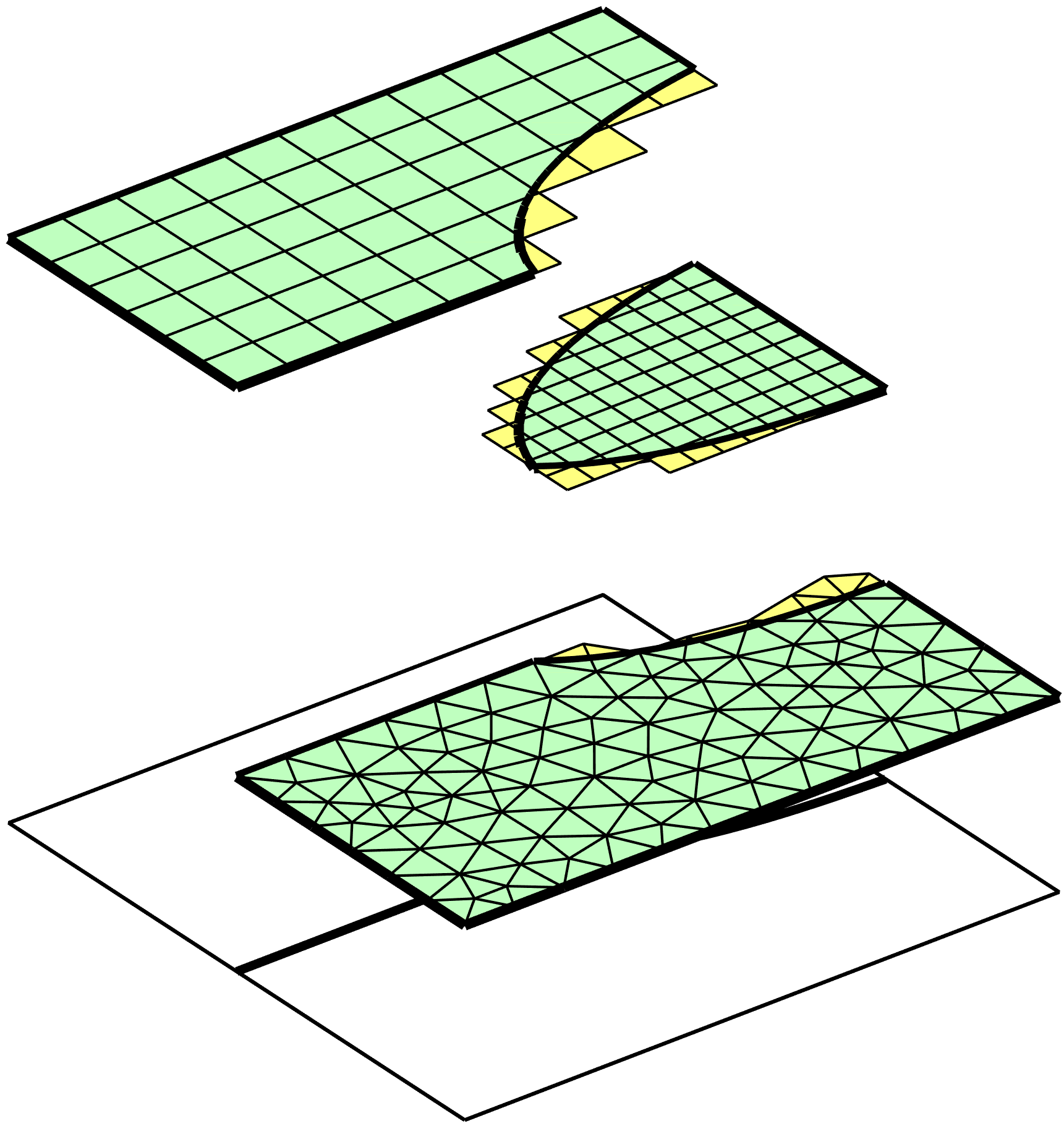}
\subcaption{Subdomain meshes}
\end{subfigure}
\begin{subfigure}[t]{0.35\linewidth}\centering
\includegraphics[width=0.9\linewidth]{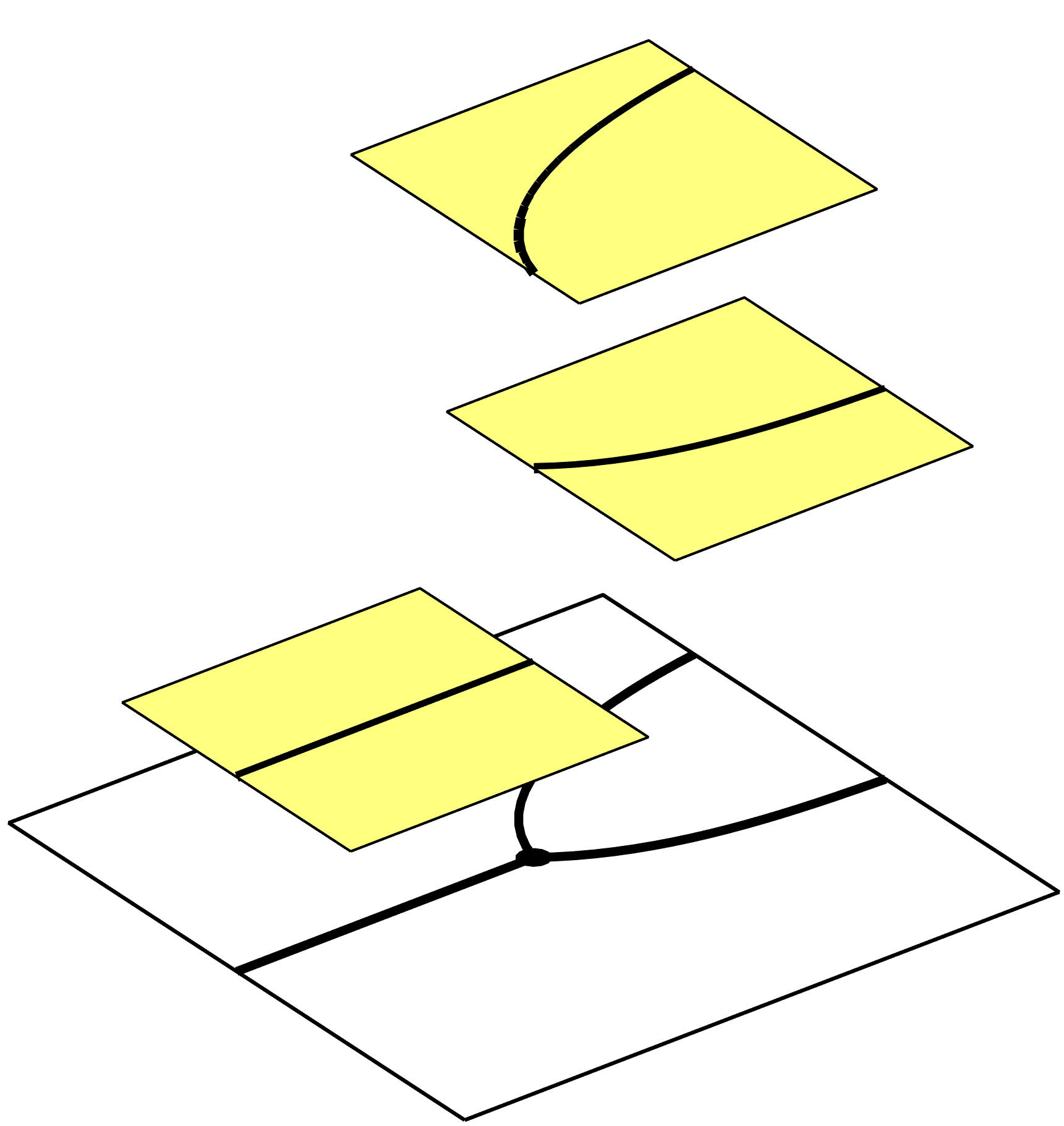}
\subcaption{Skeleton subdomain meshes}
\end{subfigure}
\caption{Meshes in the three subdomains example with skeleton subdomains embedded in a single element. (a) Subdomain solutions are approximated using $Q_2$ and $P_2$ elements on quadrilateral respectively triangle meshes. (b) Each skeleton subdomain is embedded in a single $Q_4$ element.}
\label{fig:three-parts-nonmatching-meshes}
\end{figure}

\begin{figure}\centering
\begin{subfigure}[t]{0.40\linewidth}\centering
\includegraphics[width=0.9\linewidth]{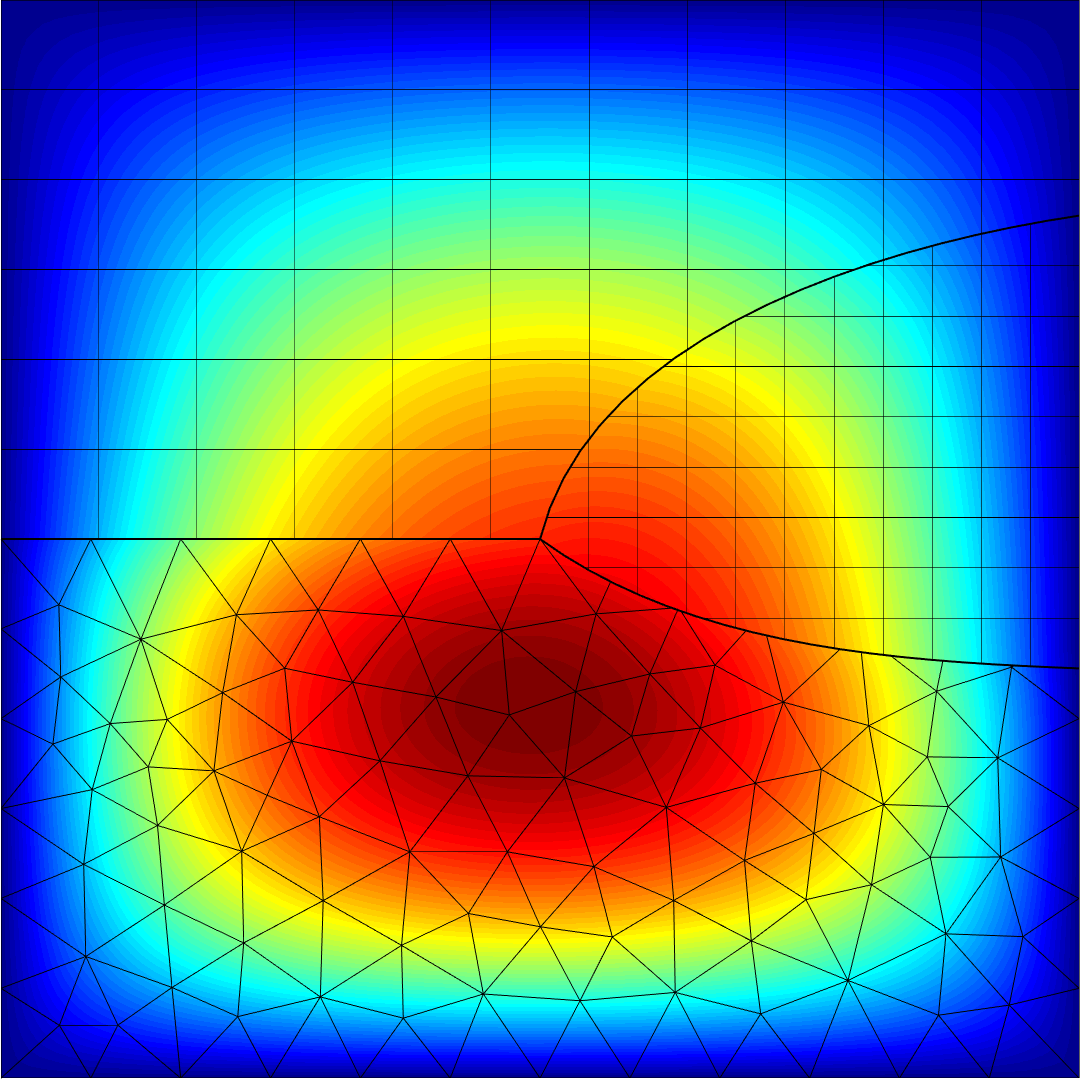}
\subcaption{Solution}
\end{subfigure}
\begin{subfigure}[t]{0.40\linewidth}\centering
\includegraphics[width=0.9\linewidth]{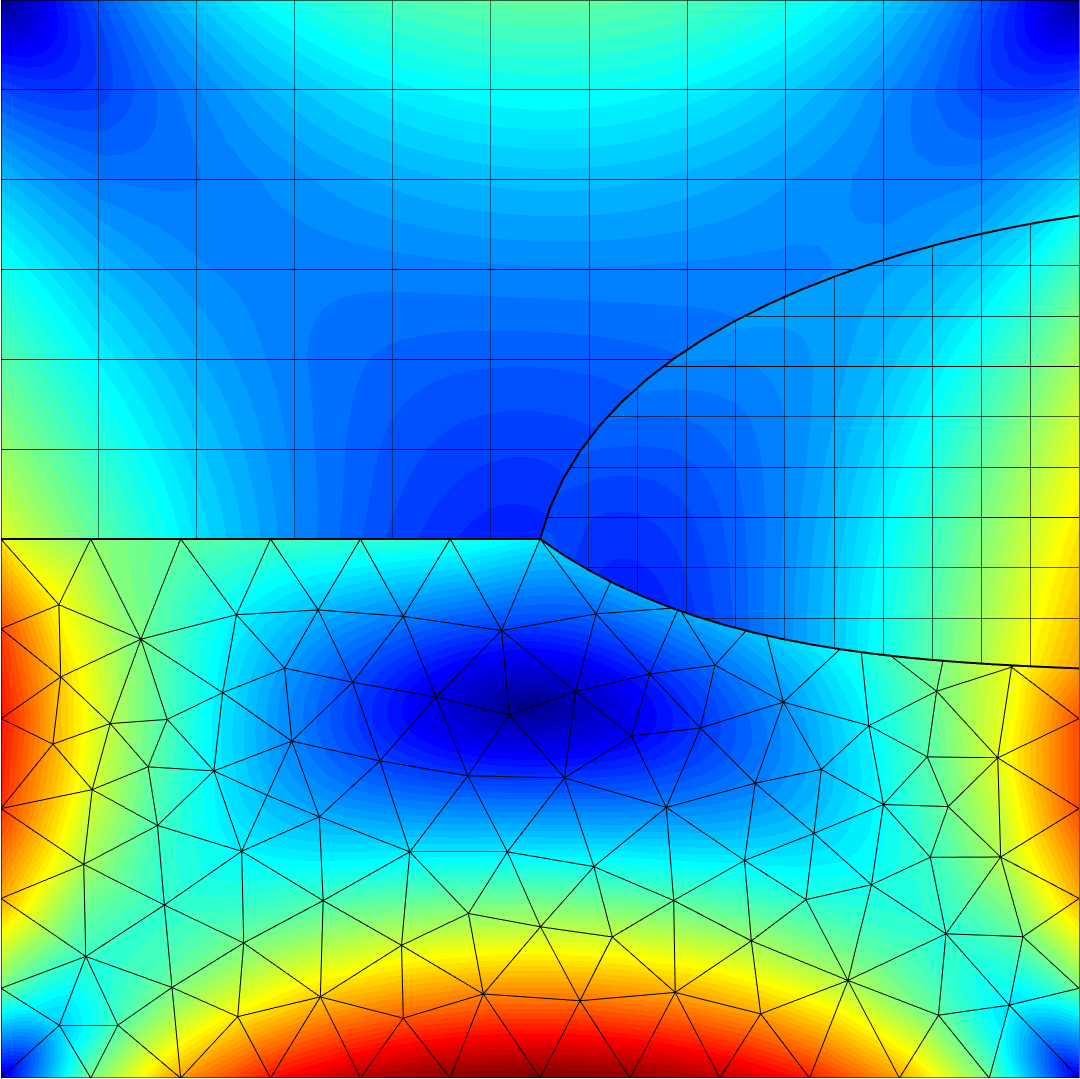}
\subcaption{Gradient magnitude}
\end{subfigure}
\caption{Three subdomains with different material coefficients and nonmatching meshes. Here meshes with $Q_2$/$P_2$ elements are used on each subdomain. On each skeleton subdomain the solution is approximated using a single $Q_4$ element.}
\label{fig:three-parts-nonmatching}
\end{figure}

For our first numerical example we partition the unit square into three subdomains with a different constant material coefficient in each subdomain, see Figure~\ref{fig:model-problem}. Thus, we in this problem have three skeleton subdomains. We consider both of the mesh constructions described in Section~\ref{section:implemention}:
\begin{itemize}
\item \emph{Global Background Grid.\ }
Here all meshes are extracted from the same background grid, see Figure~\ref{fig:three-parts-matching-meshes}, and we use $Q_2$ elements on each mesh. Note that all subdomains have cut elements and that some skeleton subdomains are curved within elements. In this setting there are no locking effects due to non-matching approximation spaces when choosing the penalty parameter $\beta$ large. A sample solution and the magnitude of its gradient are presented in Figure~\ref{fig:three-parts-matching}.

\item \emph{Single Element Interfaces.\ }
Here the mesh on each subdomain is constructed independently, some as quadrilateral meshes and some as triangular, and we equip all subdomain meshes with $Q_2$/$P_2$ elements. On each skeleton subdomain we use a single $Q_4$ element. Sample meshes in this set-up are visualized in Figure~\ref{fig:three-parts-nonmatching-meshes} and the corresponding numerical solution is presented in Figure~\ref{fig:three-parts-nonmatching}.

\end{itemize}

\begin{figure}\centering
\begin{subfigure}[t]{0.45\linewidth}\centering
\includegraphics[width=0.9\linewidth]{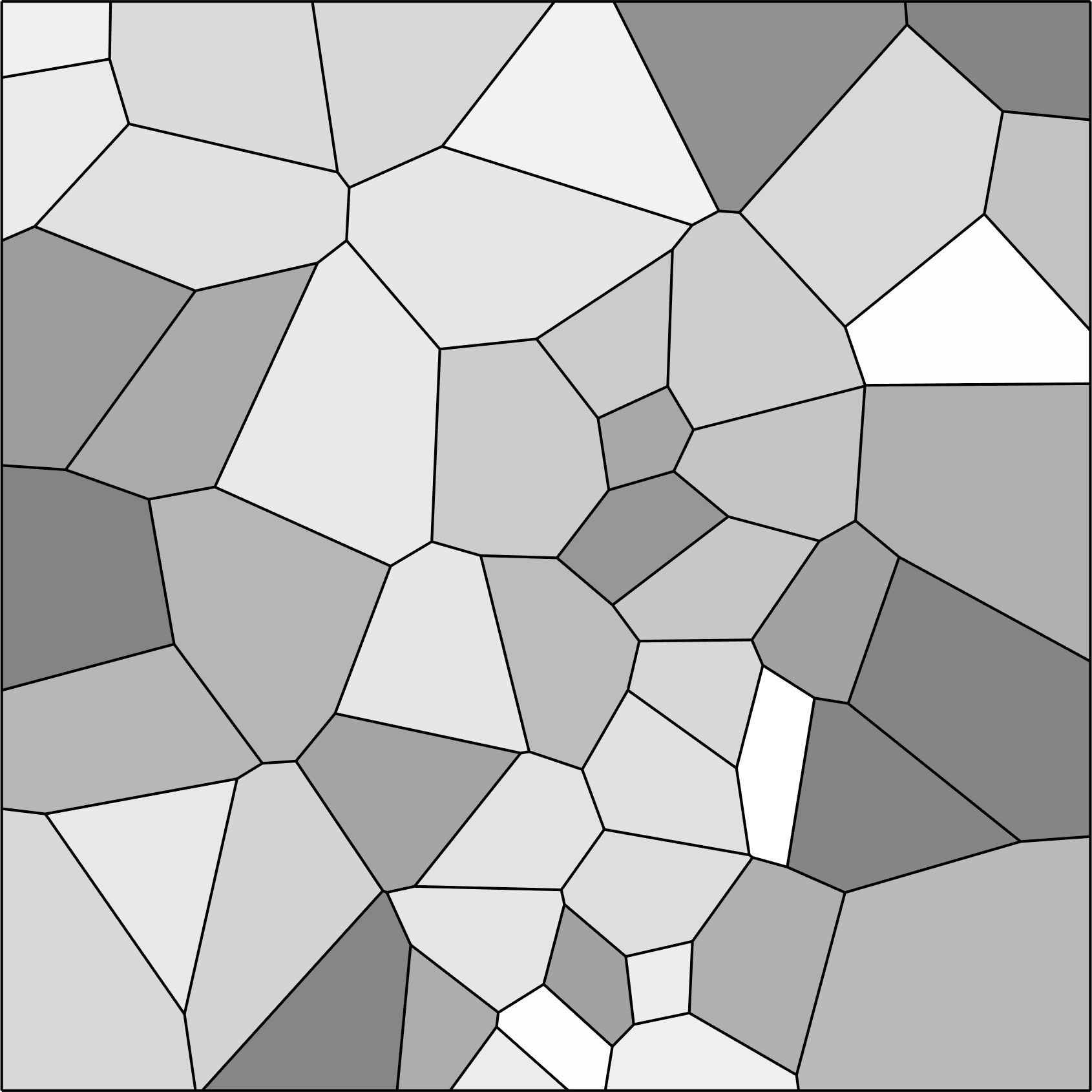}
\subcaption{Subdomain-wise coefficient $a$} \label{fig:subdomain-wise}
\label{fig:voronoi-mesh-a}
\end{subfigure}
\begin{subfigure}[t]{0.45\linewidth}\centering
\includegraphics[width=0.9\linewidth]{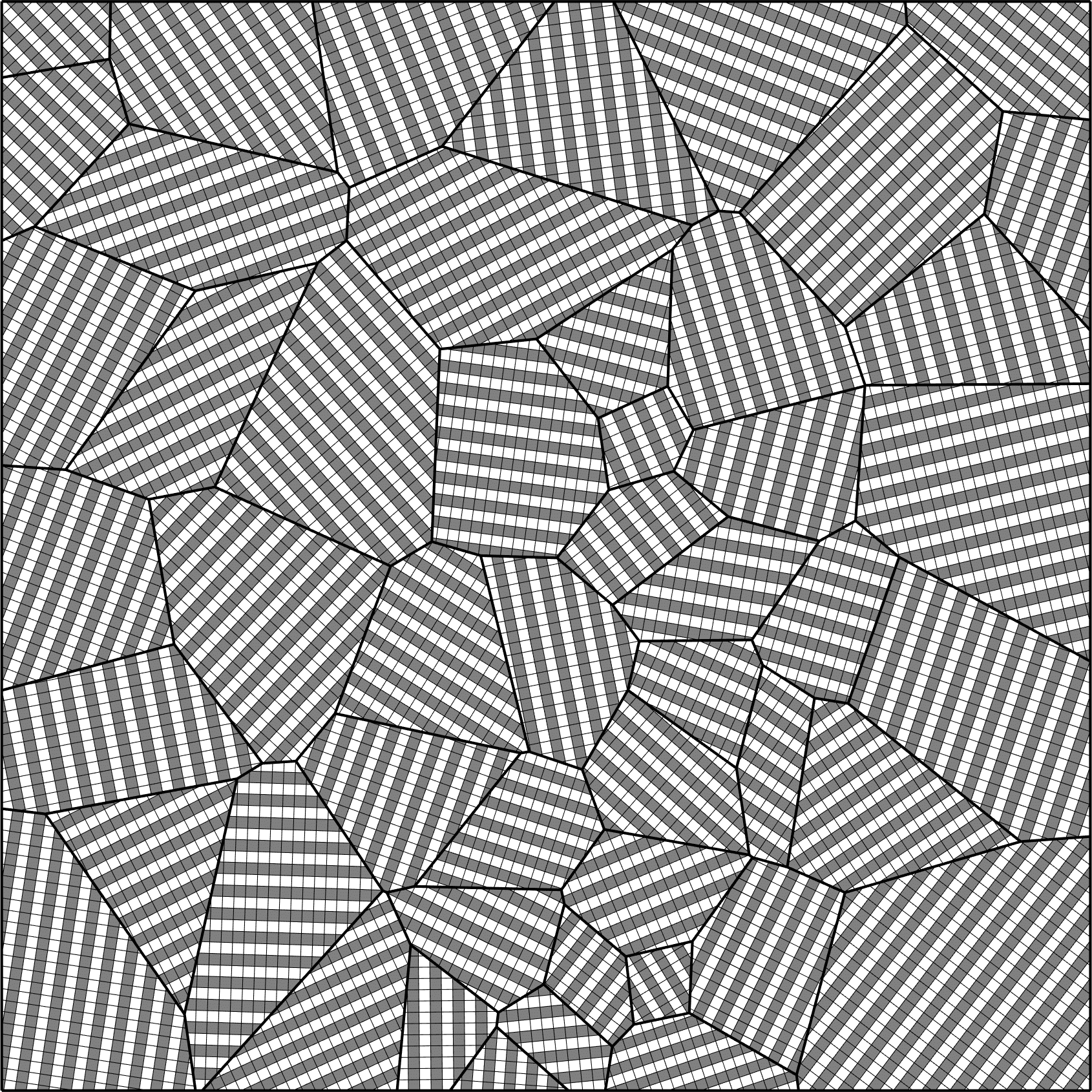}
\subcaption{Element-wise coefficient $a$} \label{fig:element-wise}
\label{fig:voronoi-mesh-b}
\end{subfigure}
\caption{Subdivisions of the unit square $[0,1]^2$ generated from Voronoi diagrams featuring varying material coefficients. (a) Domain with material coefficient $a \in [0.01,1]$ which is constant within each subdomain and chosen using a uniformly distributed random variable. (b) Domain with a randomly oriented mesh in each subdomain and a material coefficient $a$ which alternates between 1 and 1000 row-wise in the mesh.}
\label{fig:voronoi-mesh}
\end{figure}

\begin{figure}\centering
\begin{subfigure}[t]{0.35\linewidth}\centering
\includegraphics[width=0.9\linewidth]{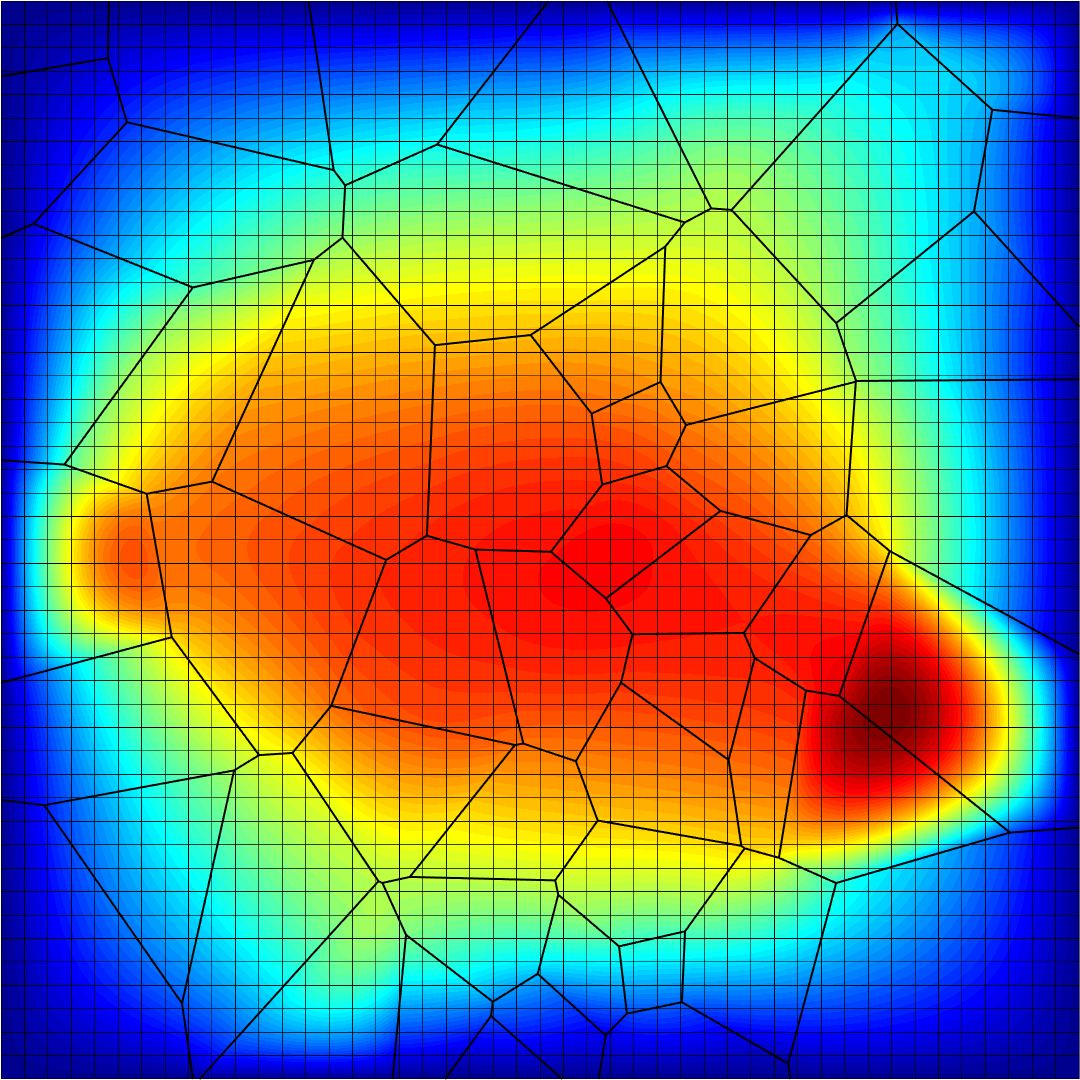}
\subcaption{$u_h$, $h=2^{-5}$}
\label{fig:variable-a}
\end{subfigure}
\begin{subfigure}[t]{0.35\linewidth}\centering
\includegraphics[width=0.9\linewidth]{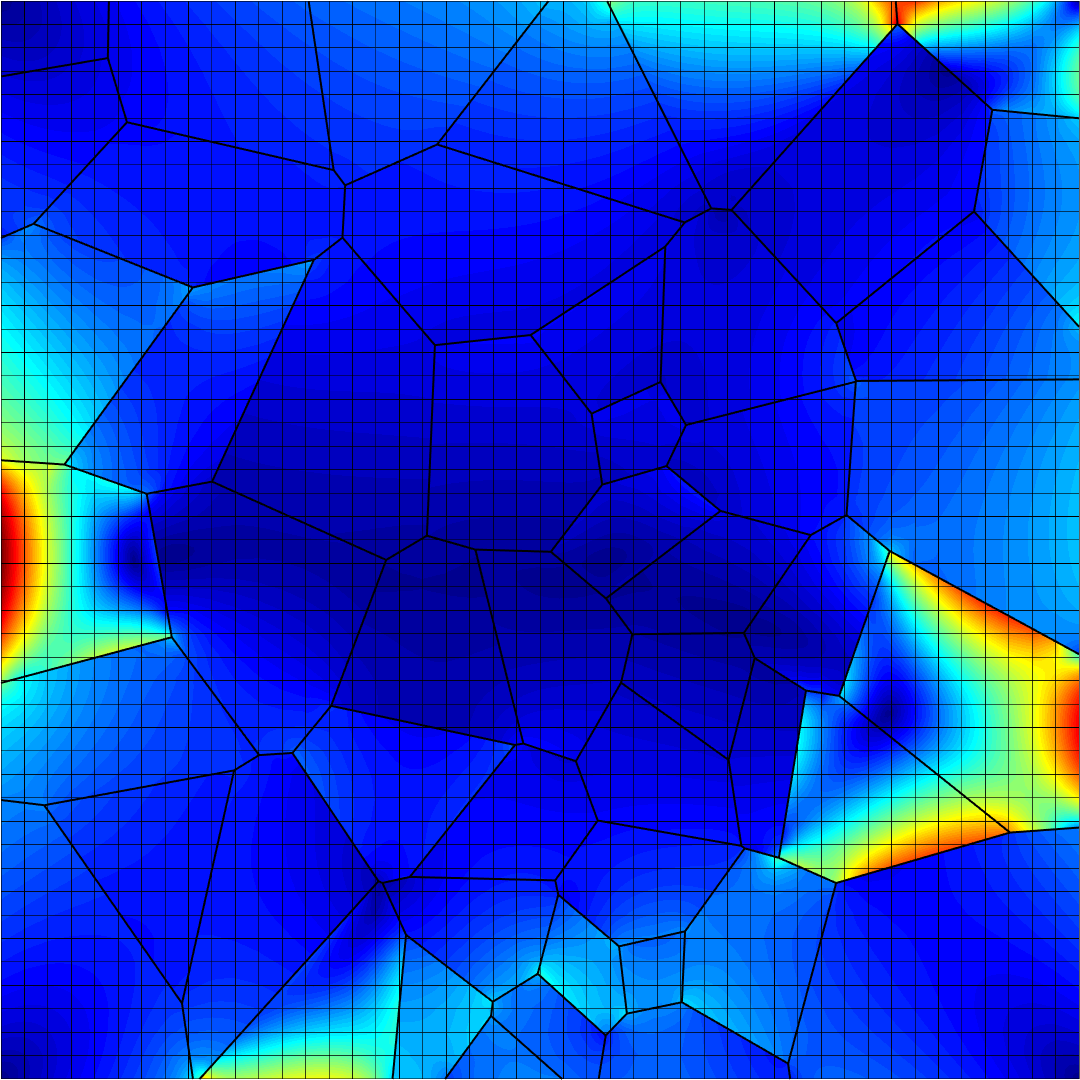}
\subcaption{$|\nabla u_h|$, $h=2^{-5}$}
\label{fig:variable-b}
\end{subfigure}

\vspace{0.5em}
\begin{subfigure}[t]{0.35\linewidth}\centering
\includegraphics[width=0.9\linewidth]{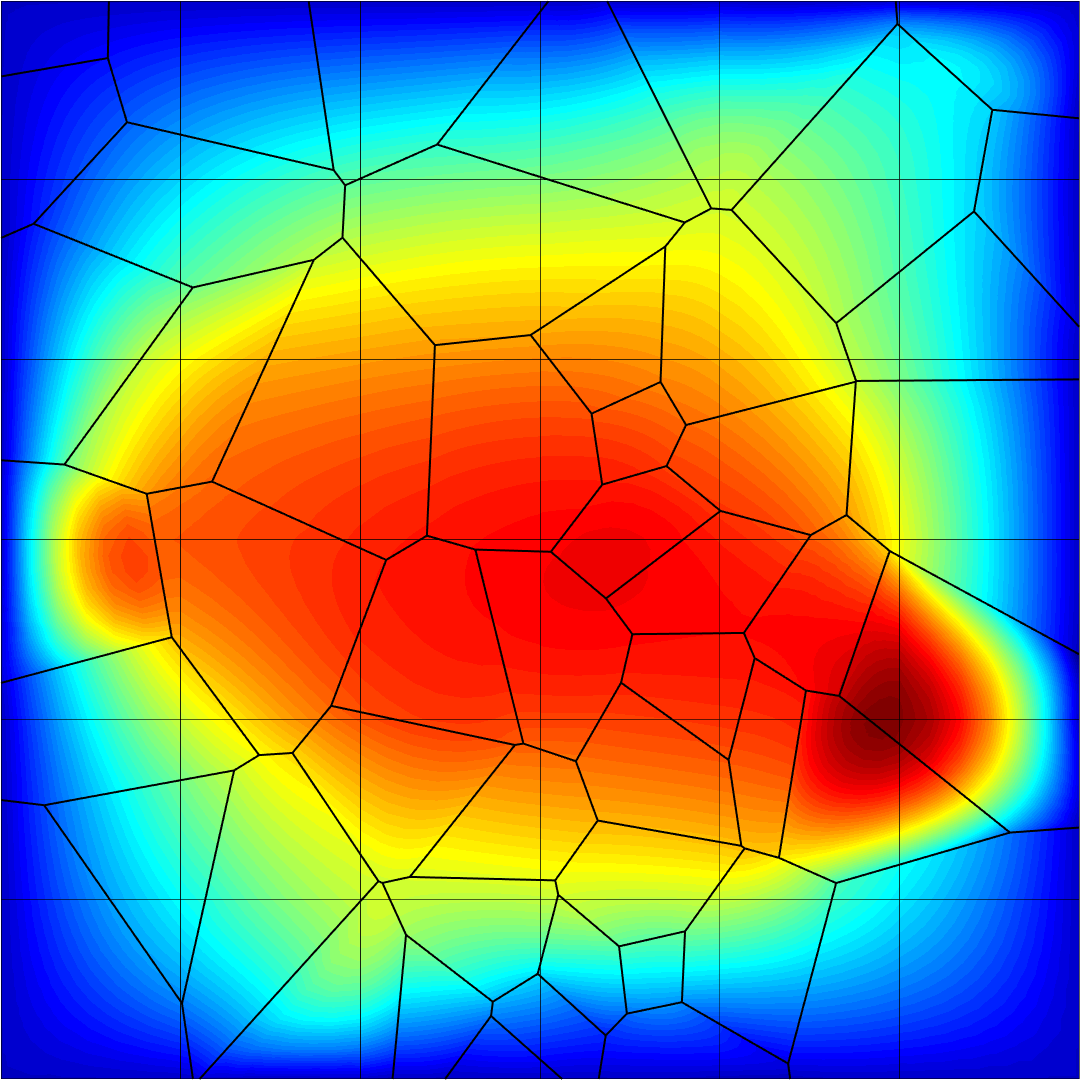}
\subcaption{$u_h$, $h=2^{-2}$}
\label{fig:variable-c}
\end{subfigure}
\begin{subfigure}[t]{0.35\linewidth}\centering
\includegraphics[width=0.9\linewidth]{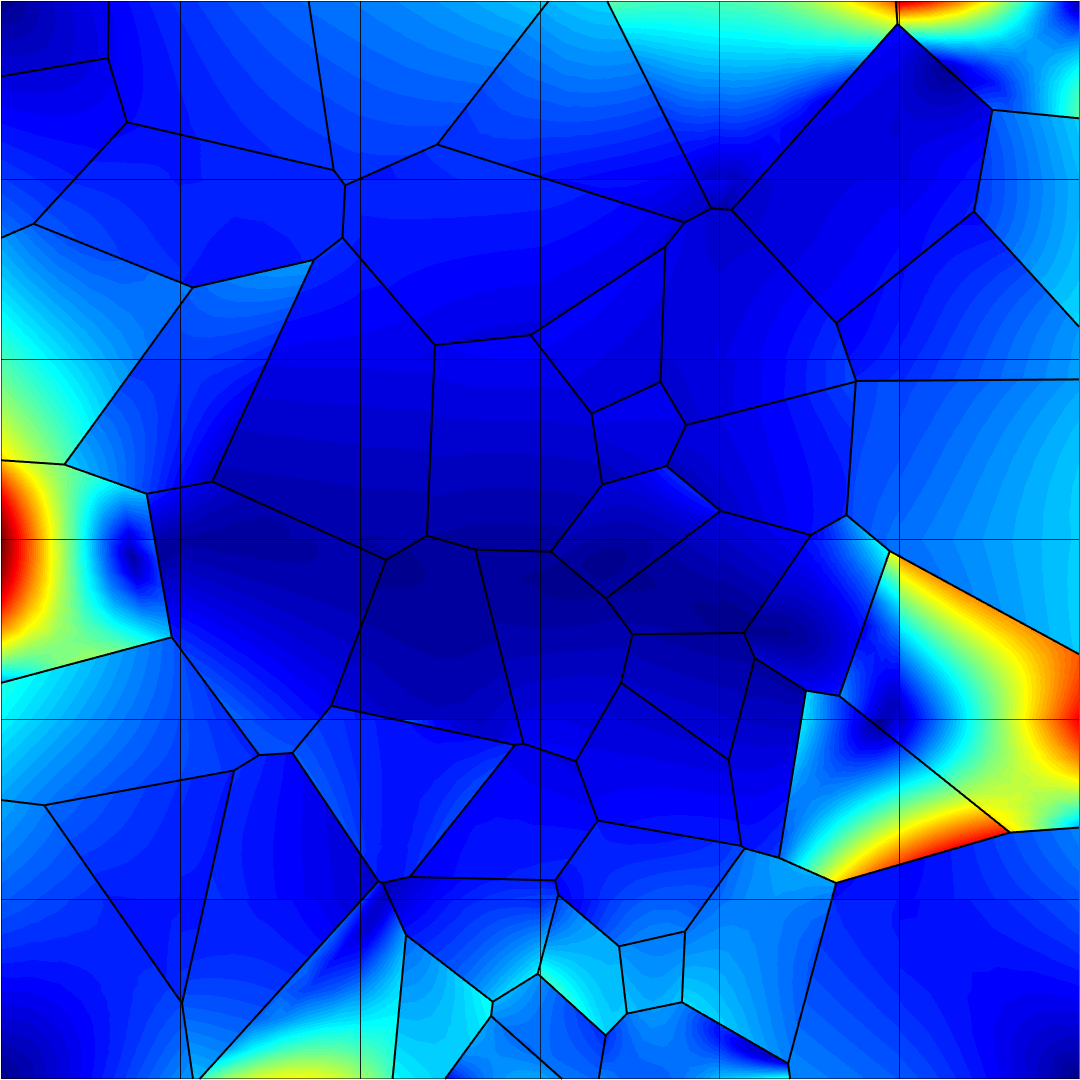}
\subcaption{$|\nabla u_h|$, $h=2^{-2}$}
\label{fig:variable-d}
\end{subfigure}

\vspace{0.5em}
\begin{subfigure}[t]{0.35\linewidth}\centering
\includegraphics[width=0.9\linewidth]{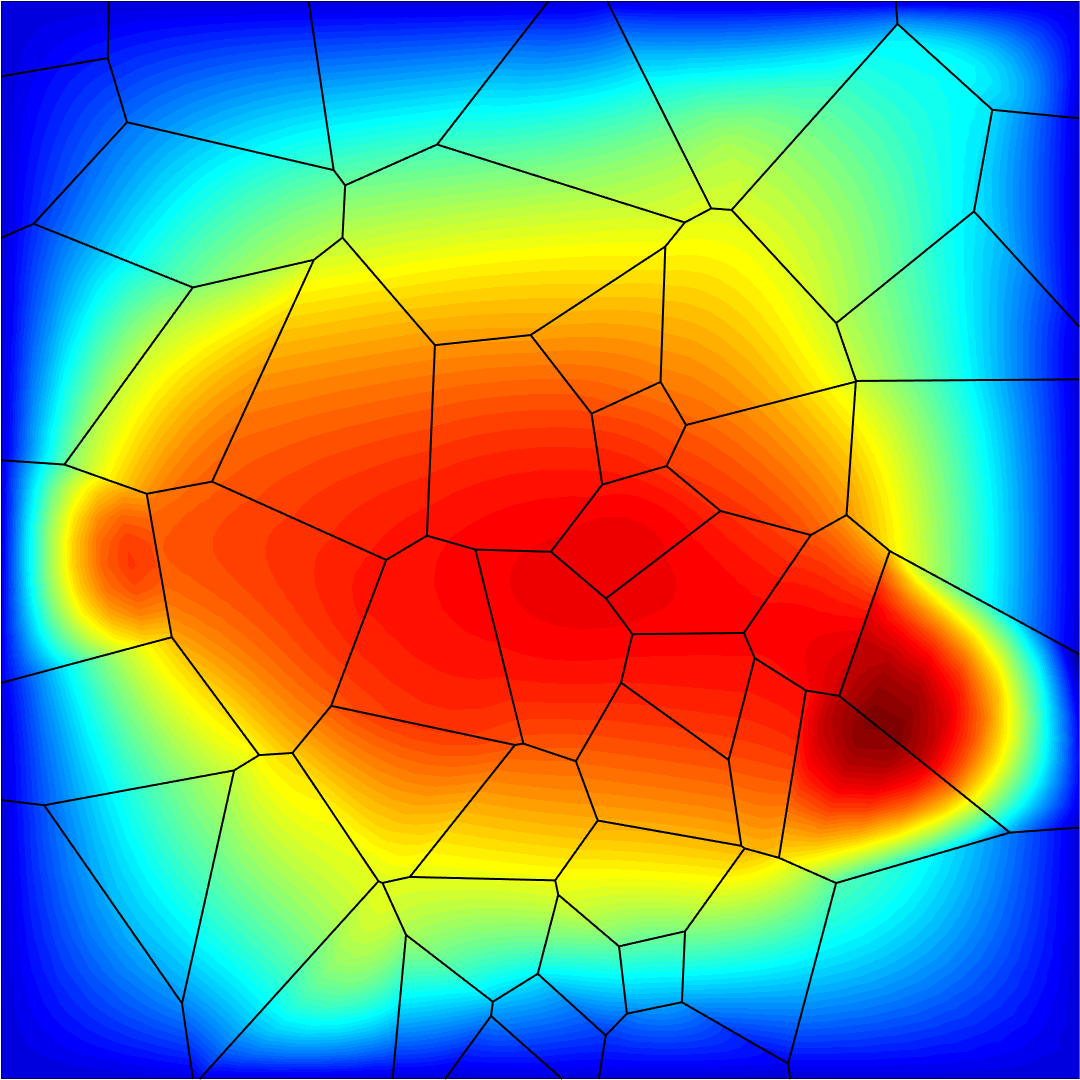}
\subcaption{$u_h$, $h=1$}
\label{fig:variable-e}
\end{subfigure}
\begin{subfigure}[t]{0.35\linewidth}\centering
\includegraphics[width=0.9\linewidth]{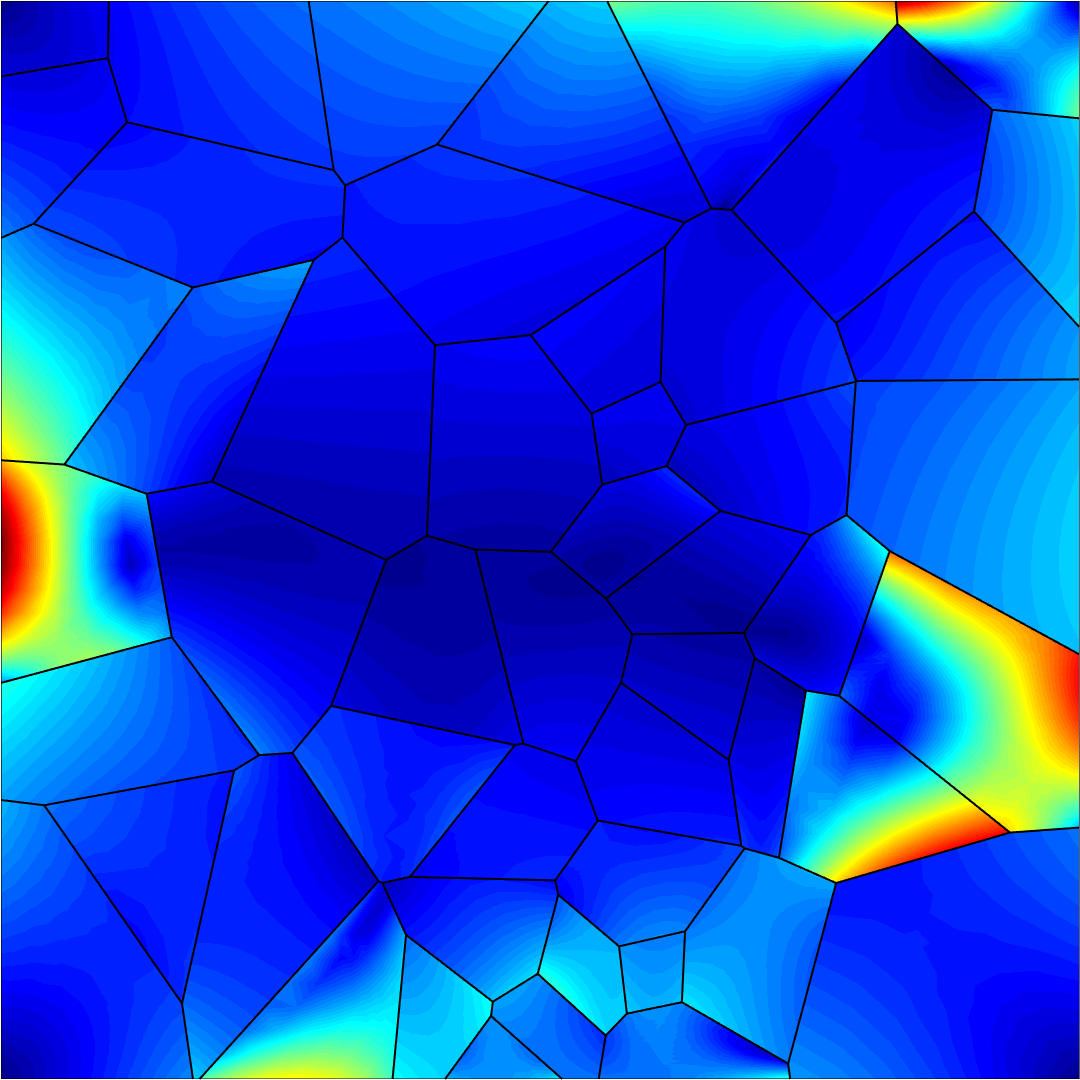}
\subcaption{$|\nabla u_h|$, $h=1$}
\label{fig:variable-f}
\end{subfigure}
\caption{Numerical solution $u_h$ and gradient magnitude $|\nabla u_h|$ on a Voronoi diagram subdivision with subdomain-wise constant material coefficient, see Figure~\ref{fig:subdomain-wise}. (a)--(b) $Q_2$ elements on meshes generated from one fine grid. (c)--(d) $Q_2$ elements on meshes generated from one coarse grid with a mesh size in the same order as the subdomain sizes. (e)--(f) A single $Q_2$ element on each subdomain and skeleton subdomain.}
\label{fig:variable}
\end{figure}

\begin{figure}\centering
\begin{subfigure}[t]{0.4\linewidth}\centering
\includegraphics[width=0.9\linewidth]{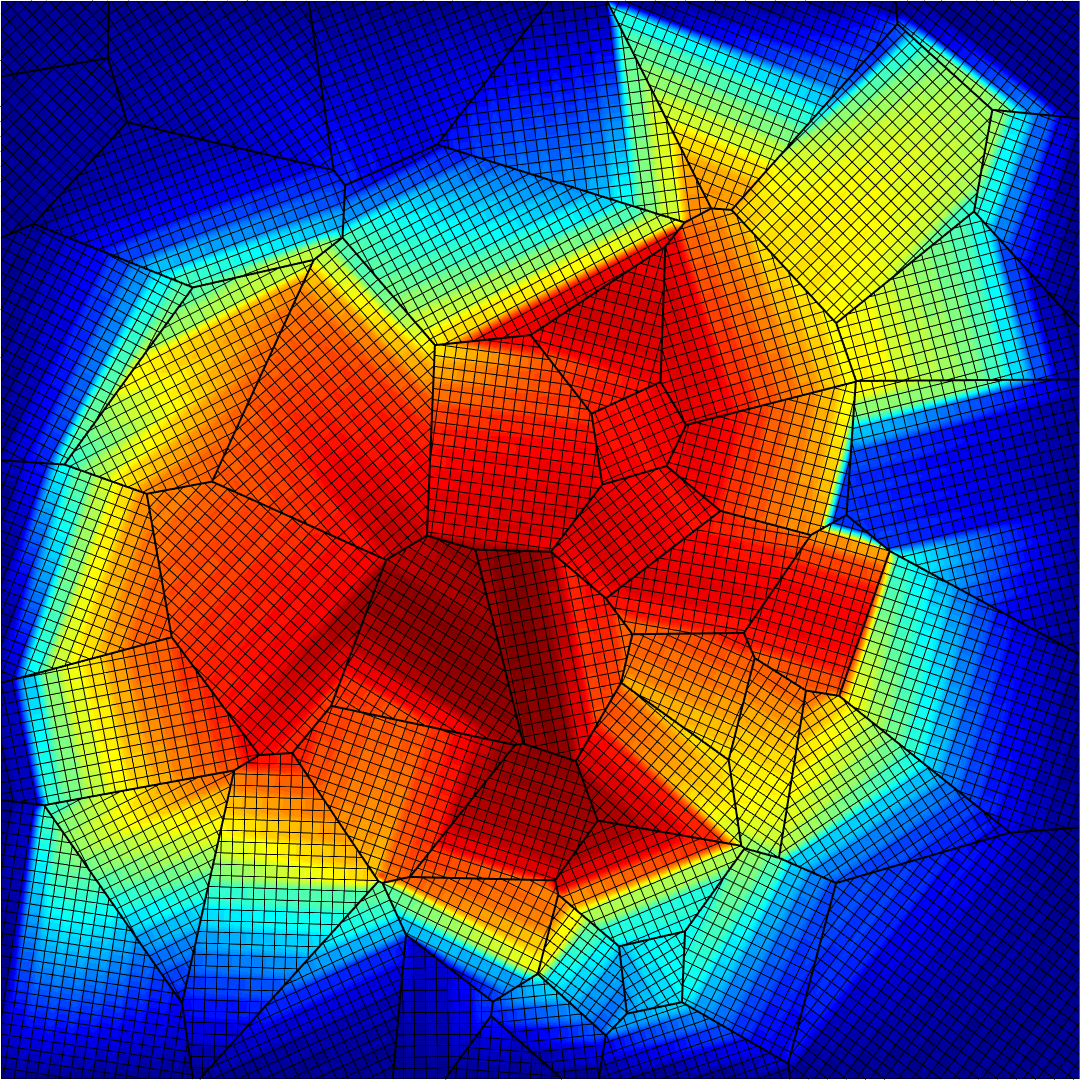}
\subcaption{Solution}
\end{subfigure}
\begin{subfigure}[t]{0.4\linewidth}\centering
\includegraphics[width=0.9\linewidth]{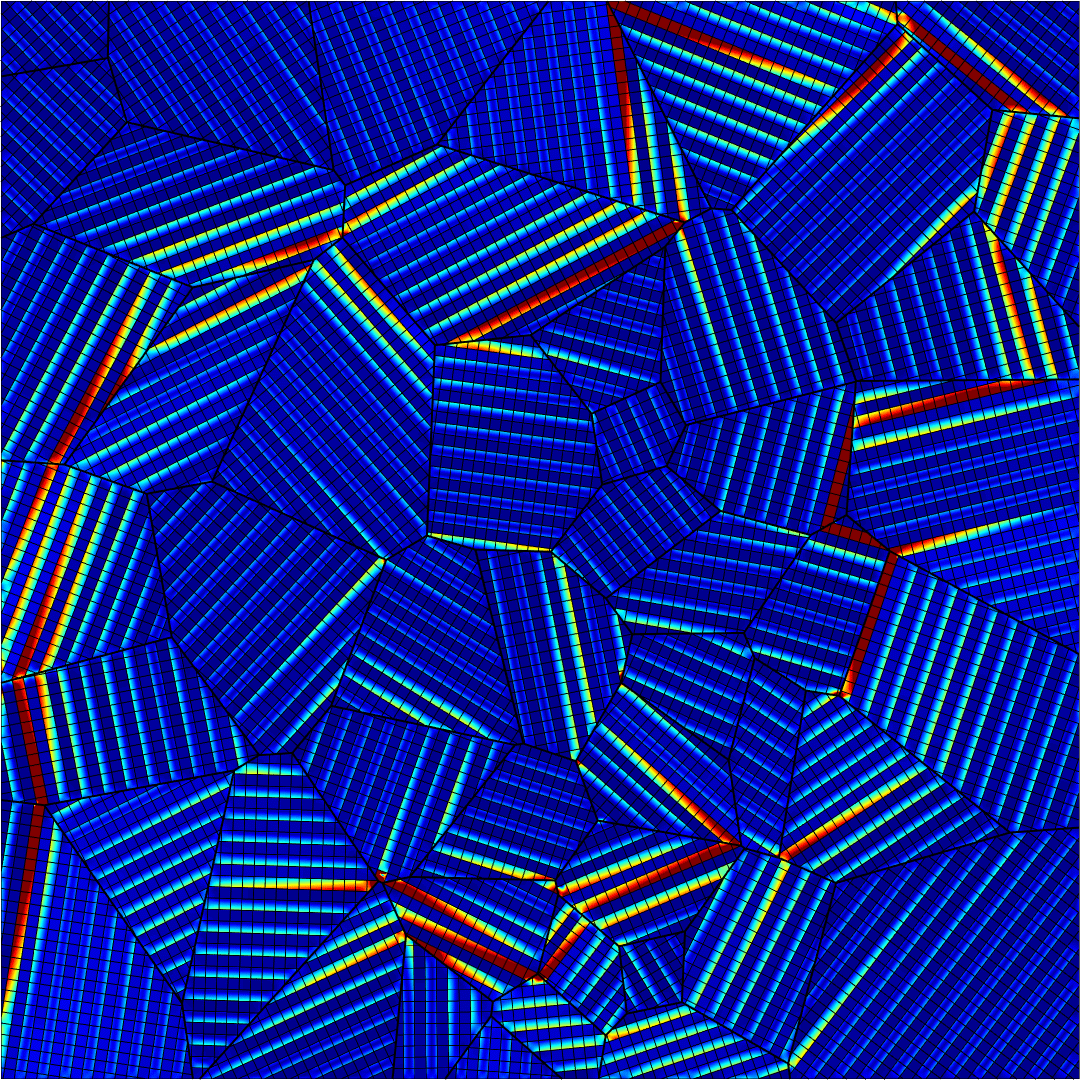}
\subcaption{Gradient magnitude}
\end{subfigure}
\caption{Numerical solution $u_h$ and gradient magnitude $|\nabla u_h|$ on a Voronoi diagram subdivision with a fine scale material coefficient pattern. On each subdomain a randomly oriented $Q_2$ mesh with a row-wise alternating material coefficient is set-up, see Figure~\ref{fig:element-wise}. The numerical solution on each skeleton subdomain is instead approximated by a single $Q_4$ element.}
\label{fig:multiscale-sol}
\end{figure}

\paragraph{Example 2: Voronoi Diagram.}
In our second example we construct a subdivision of the unit square by generating a Voronoi diagram from 50 uniformly distributed random points in $[0,1]^2$ and taking the restriction of this diagram to the unit square, see Figure~\ref{fig:voronoi-mesh}. We again consider the two different set-ups regarding mesh construction albeit we make different choices for the material coefficient in the two cases.
\begin{itemize}
\item \emph{Global Background Grid.\ }
Here we extract all meshes from the same background grid and we equip our meshes with Lagrange $Q_2$ elements. The material coefficient is constant on each subdomain $\Omega_i$ and is chosen as
$a_i = 0.01 + X$ where $X \in [0,1]$ is a uniformly distributed random variable.
This set-up is illustrated in Figure~\ref{fig:voronoi-mesh-a}.
We can easily generate background grids of any mesh size and in Figure~\ref{fig:variable} we present results for three different mesh sizes. In Figures~\ref{fig:variable-c}--\ref{fig:variable-d} we note that the method also works well when the mesh size is of the same order as the subdomain sizes. The extreme case where we only have a single element on each subdomain is presented in Figures~\ref{fig:variable-e}--\ref{fig:variable-f}. This is much like a hybridized version of so called polygonal/polyhedral elements, see the overview in \cite{polygonal2014}.
Note that in this extreme case we construct the single elements such that they are as small as possible while still containing its subdomain to avoid conditioning problems.

\item \emph{Single Element Interfaces.\ }
In the situation where we equip skeleton subdomain with a single Lagrange $Q_4$ element we choose another set-up regarding meshes and material coefficient in each bulk subdomain. We randomly orient a fine mesh equipped with Lagrange $Q_2$ elements on each subdomain and we let the material coefficient alternate between 1 and 1000 row-wise in the mesh. This set-up is illustrated in Figure~\ref{fig:voronoi-mesh-b} and the numerical solution is presented in Figure~\ref{fig:multiscale-sol}.
Note that this case is mainly an illustration of how we conveniently
can implement cases where the subdomains are defined via mappings, in this case a rotation. Of course we here loose fine scale information across the skeleton and we have made no special adaption to handle the large variation in the material coefficient.
\end{itemize}

\subsection{Convergence Studies}

\begin{figure}\centering
\begin{subfigure}[t]{0.3\linewidth}\centering
\includegraphics[width=0.9\linewidth]{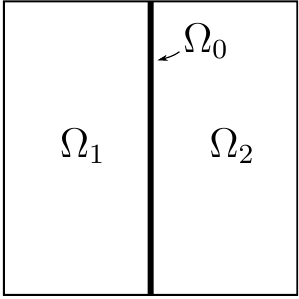}
\subcaption{Domain}
\end{subfigure}
\begin{subfigure}[t]{0.3\linewidth}\centering
\includegraphics[width=0.9\linewidth]{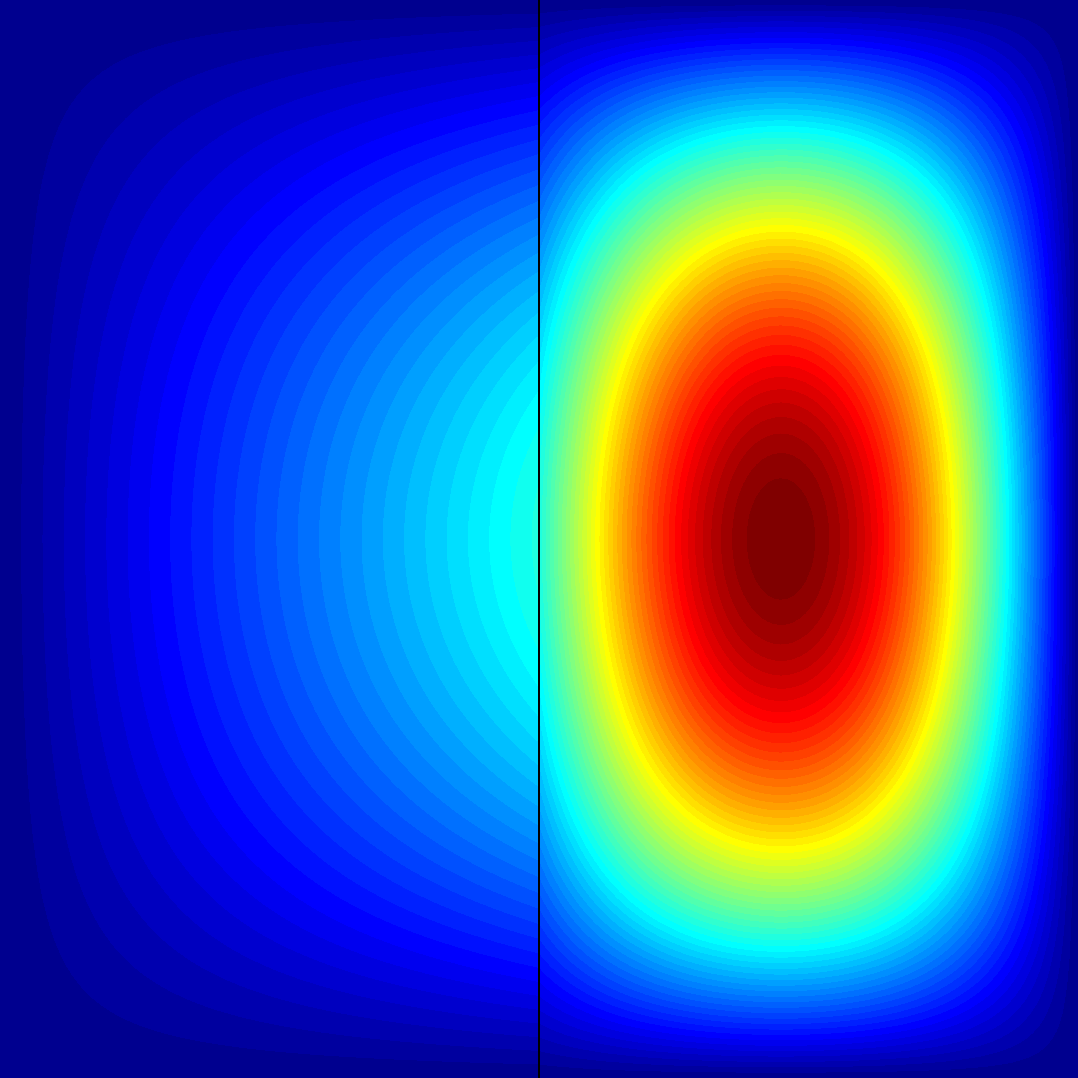}
\subcaption{Solution}
\end{subfigure}
\begin{subfigure}[t]{0.3\linewidth}\centering
\includegraphics[width=0.9\linewidth]{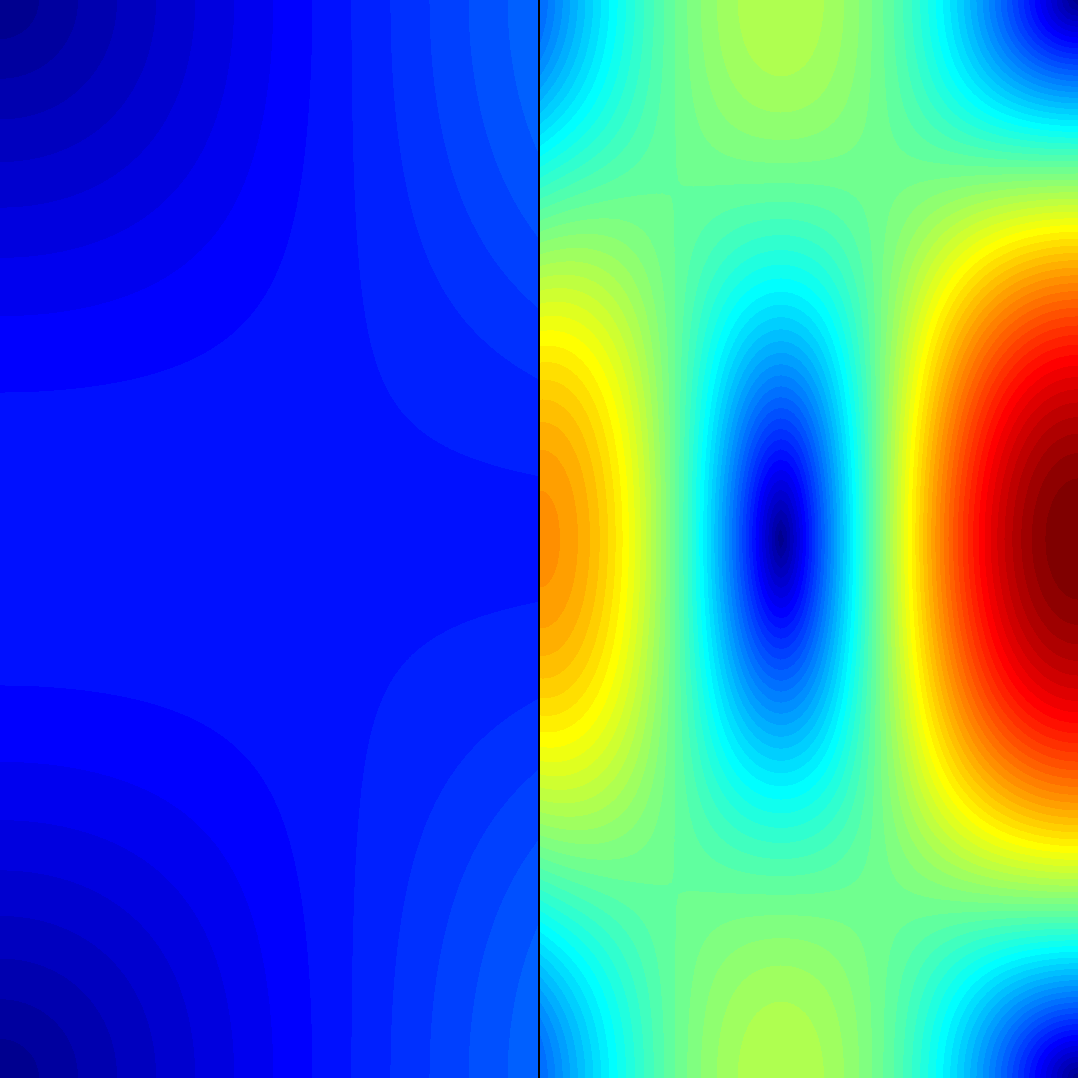}
\subcaption{Gradient magnitude}
\end{subfigure}
\caption{Problem with known exact solution used in convergence studies. (a) The domain is the unit square $[0,1]^2$ partitioned into two mirror symmetric subdomains according to the figure with material coefficients $a_1=1$ and $a_2=2\pi-1$. (b) The exact solution \eqref{eq:exact-solution}. (c) The gradient magnitude of the exact solution.}
\label{fig:convergence-setup}
\end{figure}

To study the convergence of the method in energy and $L2$ norms we partition the unit square $\Omega$ into its left and right halves, $\Omega_1$ and $\Omega_2$, and we manufacture a problem with known exact solution from the ansatz
\begin{align} \label{eq:exact-solution}
u =
\left\{\,
u_0 = \frac{1}{2}\sin(\pi y)
\,,\
u_1 = x\sin(\pi y)
\,,\
u_2 = (1 - x - \sin(2\pi x))\sin(\pi y)
\,\right\}
\end{align}
with coefficients $a_1 = (2\pi -1),\ a_2 = 1$. The domain, exact solution and exact gradient magnitude for this problem are displayed in Figure~\ref{fig:convergence-setup}.

\begin{figure}\centering
\begin{subfigure}[t]{0.35\linewidth}\centering
\includegraphics[width=0.9\linewidth]{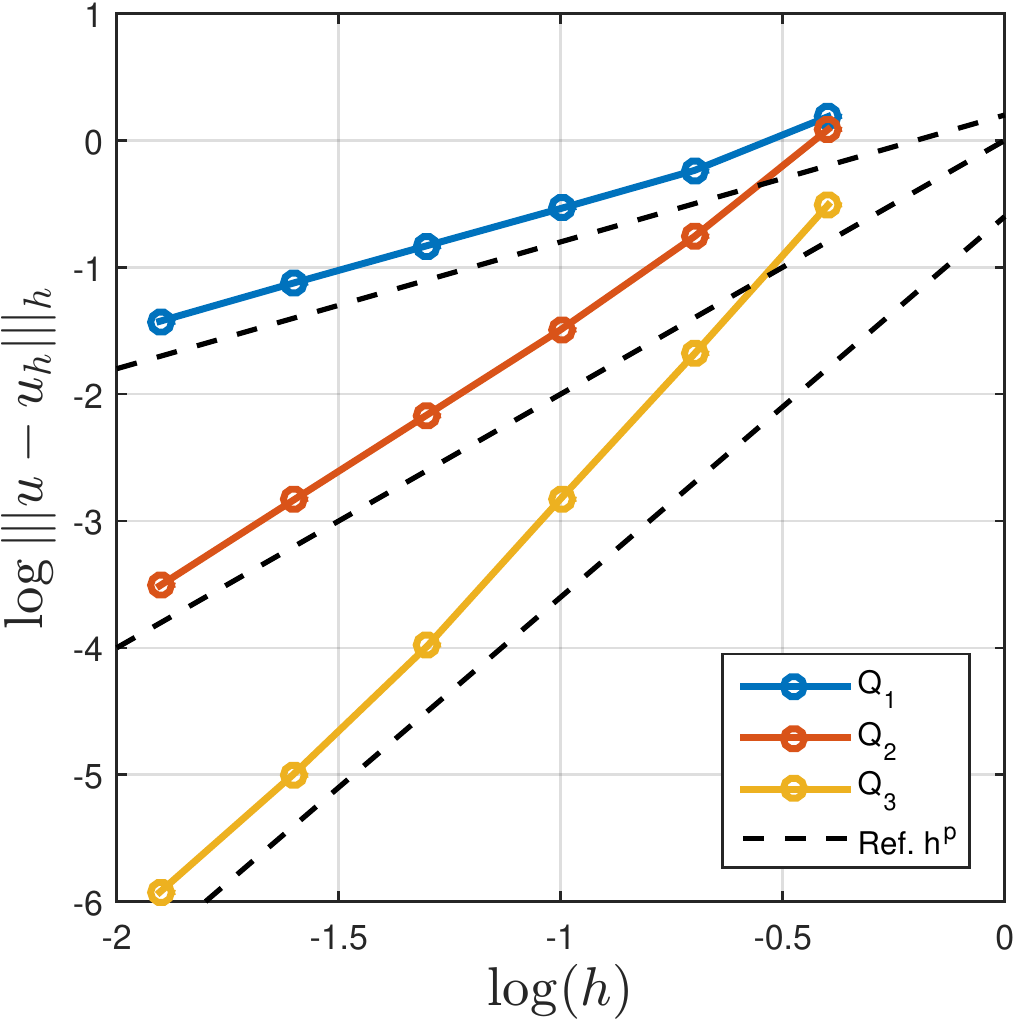}
\subcaption{Energy norm}
\end{subfigure}
\begin{subfigure}[t]{0.35\linewidth}\centering
\includegraphics[width=0.9\linewidth]{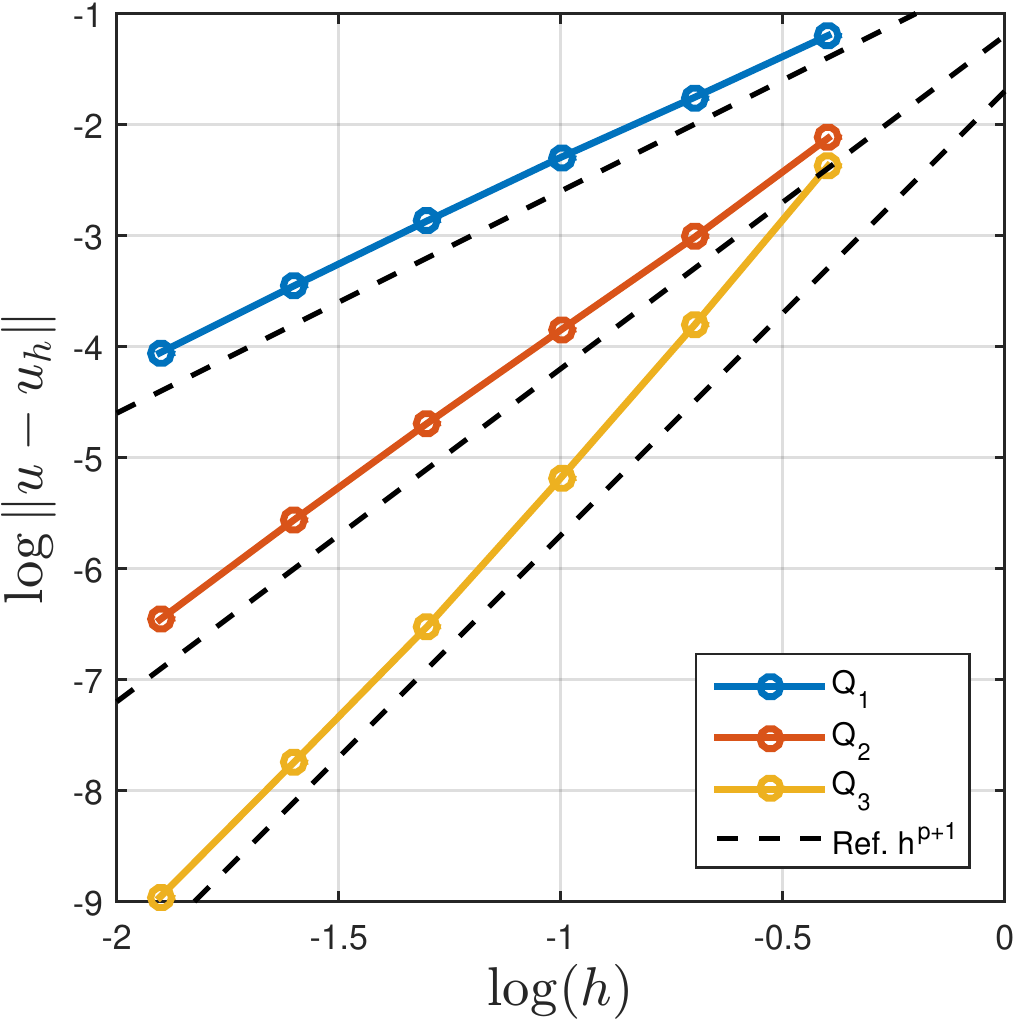}
\subcaption{$L^2$ norm}
\end{subfigure}
\caption{Convergence studies using meshes all from the same background grid. In all meshes the same elements are used ($Q_1$--$Q_3$) and we achieve what we expect to be optimal convergence rates of $O(h^{p})$ in the energy norm (a) respectively $O(h^{p+1})$ in the $L^2$ norm (b).}
\label{fig:convergence-matching}
\end{figure}

\begin{figure}\centering
\begin{subfigure}[t]{0.35\linewidth}\centering
\includegraphics[width=0.9\linewidth]{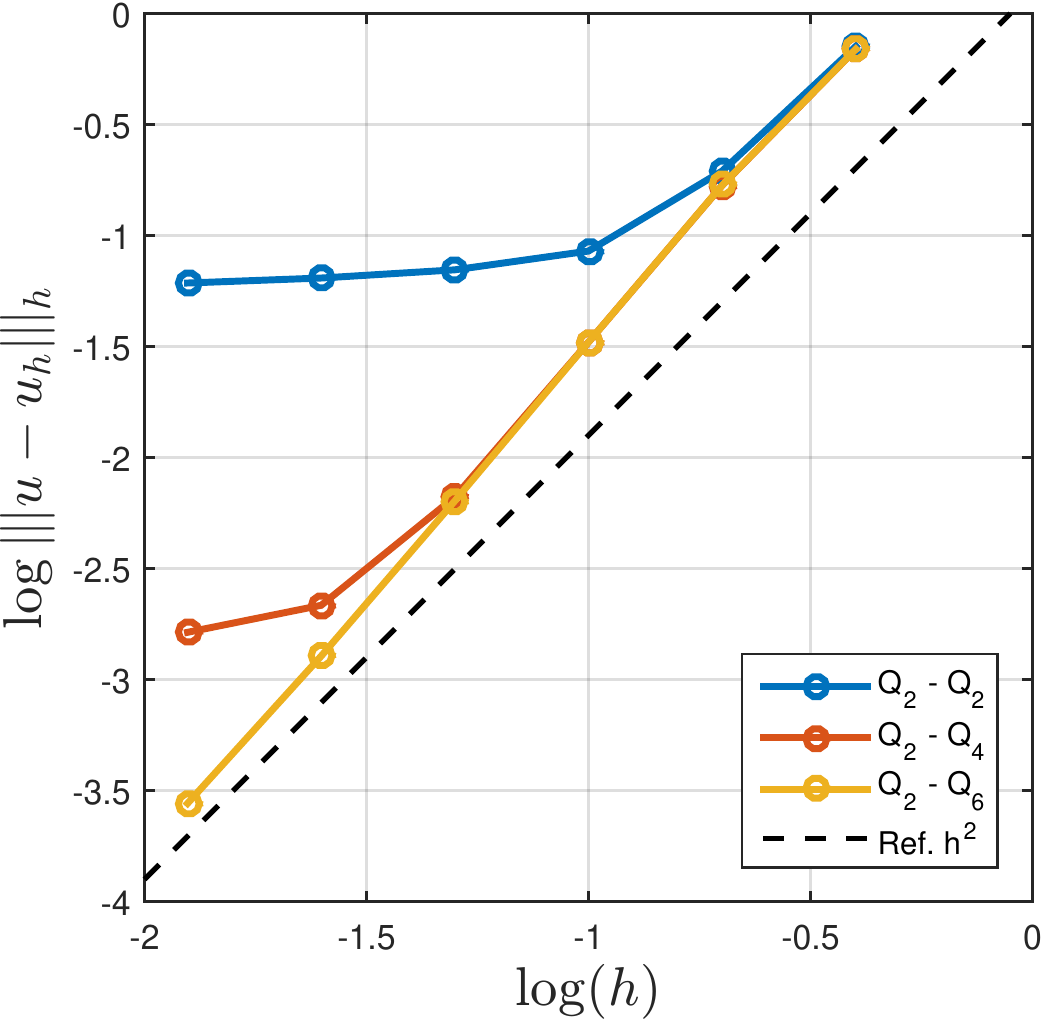}
\subcaption{Energy norm}
\end{subfigure}
\begin{subfigure}[t]{0.35\linewidth}\centering
\includegraphics[width=0.9\linewidth]{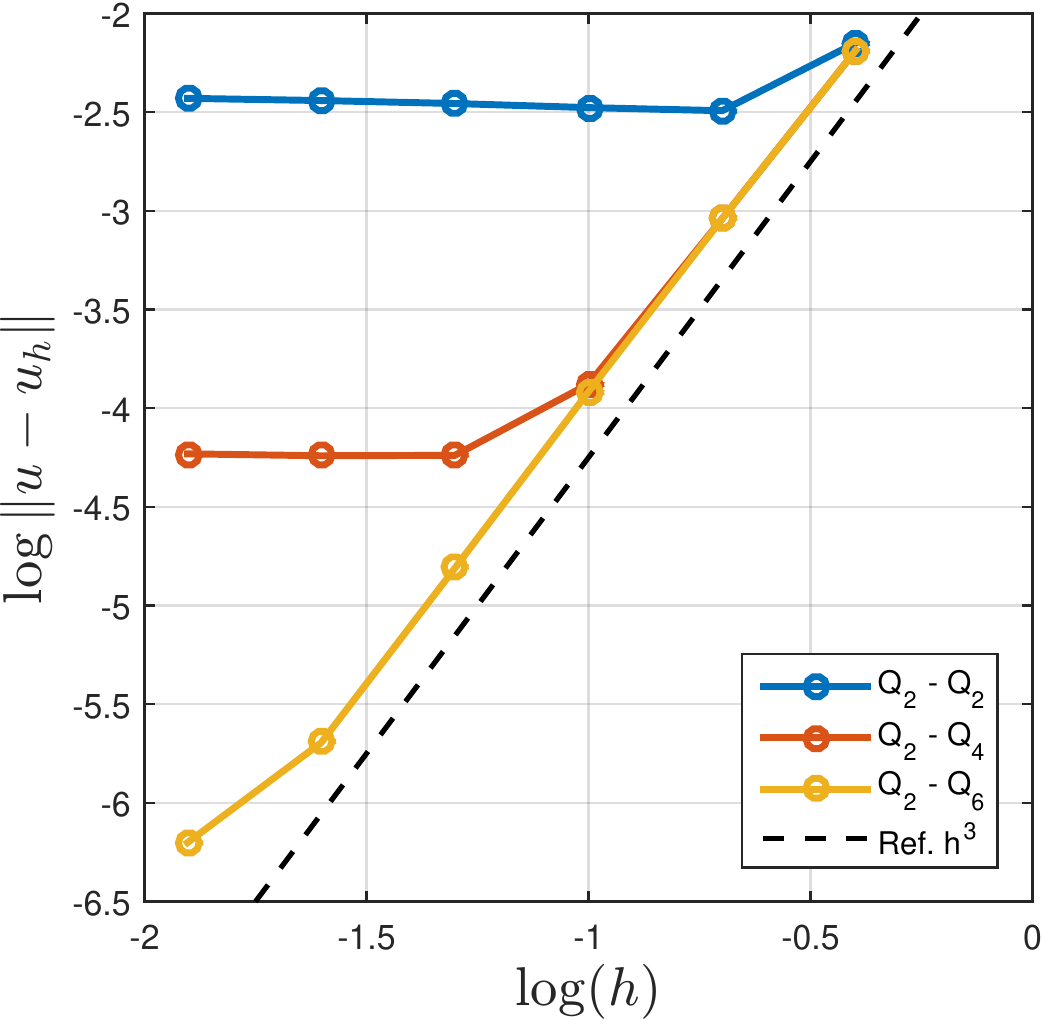}
\subcaption{$L^2$ norm}
\end{subfigure}
\caption{Convergence studies using non-matching meshes for the subdomains and a single polynomial for each skeleton subdomain. On the subdomains $Q_2$ elements are used and on the skeleton subdomains $Q_2$--$Q_6$ polynomials are used. Because there is no $h$ refinement of the skeleton subdomains the convergence levels out as the polynomial approximations on the skeleton become the dominant source of error.}
\label{fig:convergence-nonmatching}
\end{figure}

\begin{itemize}
\item \emph{Global Background Grid.\ }
In Figure~\ref{fig:convergence-matching} we present convergence results in the case where all meshes are extracted form the same quadrilateral background grid using $Q_p$ elements, $p=1,2,3$, and we achieve the expected convergence rates of $O(h^{p})$ and $O(h^{p+1})$ in the energy and the $L^2$ norm, respectively.

\item \emph{Single Element Interfaces.\ }
In this case we instead use $Q_2$ elements of size $h$ on each subdomain but only a single Lagrange $Q_p$ element, $p=2,4,6$, on each skeleton subdomain. As the skeleton subdomain meshes are not refined with smaller $h$ this naturally imposes a lower bound on the errors, which we also note in the convergence results presented in Figure~\ref{fig:convergence-nonmatching}.

\end{itemize}

\bigskip
\paragraph{Acknowledgements.}
This research was supported in part by the Swedish Foundation
for Strategic Research Grant No.\ AM13-0029, the Swedish Research
Council Grants Nos.\  2013-4708, 2017-03911, and the Swedish
Research Programme Essence. EB was supported by EPSRC research grants EP/P01576X/1 and EP/P012434/1.


\bibliographystyle{abbrv}
\footnotesize{
\bibliography{hybrid-refs}
}

\bigskip
\bigskip
\noindent
\footnotesize {\bf Authors' addresses:}

\smallskip
\noindent
Erik Burman,  \quad \hfill \addressuclshort\\
{\tt e.burman@ucl.ac.uk}

\smallskip
\noindent
Daniel Elfverson,  \quad \hfill \addressumushort\\
{\tt daniel.elfverson@umu.se}

\smallskip
\noindent
Peter Hansbo,  \quad \hfill \addressjushort\\
{\tt peter.hansbo@ju.se}

\smallskip
\noindent
Mats G. Larson,  \quad \hfill \addressumushort\\
{\tt mats.larson@umu.se}

\smallskip
\noindent
Karl Larsson, \quad \hfill \addressumushort\\
{\tt karl.larsson@umu.se}

\end{document}